\documentclass[11pt]{article}
\usepackage[utf8]{inputenc}
\usepackage{setspace}
\usepackage{comment}

\usepackage[left=1in,right=1in,top=1in,bottom=1in,bindingoffset=0cm]{geometry}
\usepackage{hyperref}
\hypersetup{
    colorlinks=true,
    linkcolor=blue,
    filecolor=magenta,      
    urlcolor=red,
    citecolor=gray,
    linktoc=page,
}
\usepackage{graphicx,wrapfig,lipsum,booktabs}
\graphicspath{ {images/} }
\usepackage{amsmath}
\usepackage{bm}
\usepackage{amsthm}
\usepackage{tikz}
\usepackage{pgfplots}
\pgfplotsset{compat = newest, width = 10cm}
\usepgfplotslibrary{fillbetween}
\usepackage{caption}
\usepackage{amssymb}
\usepackage{mathtools}

\usepackage{bbm}
\usepackage{enumitem}
\usepackage{algorithm}
\usepackage{algpseudocode}
\counterwithin{algorithm}{section}
\allowdisplaybreaks
\usepackage{breakurl}
\usepackage[backend=biber, style=alphabetic, sorting=nyt, maxnames=100,backref=true, firstinits=true, sortcites = true]{biblatex}
\DeclareFieldFormat[unpublished]{howpublished}{#1, }

\def\ER{Erd\H{o}s--R\'{e}nyi }

\def\E{\mathbb{E}}

\def\N{\mathbb{N}}

\def\E{\mathbb{E}}

\def\P{\mathbb{P}}

\def\eps{\varepsilon}

\def\1{\mathbf{1}}

\def\tce{t_c + \eps}
\def\tce2{t_c + \frac{\eps}{2}}

\def\Hrrnp{\mathcal{H}(r,n,p)}

\newtheorem{Theorem}{Theorem}[section]

\newtheorem{Claim}[Theorem]{Claim}
\newtheorem{Fact}[Theorem]{Fact}

\newtheorem{lemma}[Theorem]{Lemma}
\newtheorem*{lemma*}{Lemma}

\newtheorem*{Theorem*}{Theorem}

\newtheorem{Definition}[Theorem]{Definition}

\newcommand{\boldgamma}{\bm{\gamma}}
\newcommand{\set}[1]{\{#1\}}
\newcommand{\dom}{\mathrm{dom}}

\newenvironment{breakablealgorithm}
  {
   \begin{center}
     \refstepcounter{algorithm}
     \hrule height.8pt depth0pt \kern2pt
     \renewcommand{\caption}[2][\relax]{
       {\raggedright\textbf{\ALG@name~\thealgorithm} ##2\par}%
       \ifx\relax##1\relax 
         \addcontentsline{loa}{algorithm}{\protect\numberline{\thealgorithm}##2}%
       \else 
         \addcontentsline{loa}{algorithm}{\protect\numberline{\thealgorithm}##1}%
       \fi
       \kern2pt\hrule\kern2pt
     }
  }{
     \kern2pt\hrule\relax
   \end{center}
  }

\title{Balanced Colorings of \ER Hypergraphs}
\author{Abhishek Dhawan\footnote{Department of Mathematics, University of Illinois Urbana-Champaign. Partially supported by NSF RTG grant DMS-1937241. Email: adhawan2@illinois.edu.}
\and 
Yuzhou Wang\footnote{School of Mathematics, Georgia Institute of Technology. Email: ywang3694@gatech.edu.}}
\date{}

\begin{document}
\sloppy

\maketitle

\begin{abstract}
    An $r$-uniform hypergraph $H = (V, E)$ is $r$-partite if there exists a partition of the vertex set into $r$ parts such that each edge contains exactly one vertex from each part. We say an independent set in such a hypergraph is balanced if it contains an equal number of vertices from each partition. The balanced chromatic number of $H$ is the minimum value $q$ such that $H$ admits a proper $q$-coloring where each color class is a balanced independent set. In this note, we determine the asymptotic behavior of the balanced chromatic number for sparse $r$-uniform $r$-partite Erd\H{o}s--R\'enyi hypergraphs. A key step in our proof is to show that any balanced colorable hypergraph of average degree $d$ admits a proper balanced coloring with $r(r-1)d + 1$ colors. This extends a result of Feige and Kogan on bipartite graphs to this more general setting.
\end{abstract}

\subsection*{Basic Notation}

For $n \in \N$, we let $[n] \coloneqq \set{1, \ldots, n}$.
For a hypergraph $H$, its vertex and edge sets are denoted $V(H)$ and $E(H)$ respectively.
If every edge in $E(H)$ contains exactly $r$ vertices, we say the hypergraph is $r$-uniform (graphs are $2$-uniform hypergraphs).
For a set $X$ and $r \in \N$, we let $\binom{X}{r}$ be the collection of $r$-element subsets of $X$.
Let $H = (V, E)$ be an $r$-uniform hypergraph, i.e., $E \subseteq \binom{V}{r}$.
We say $H$ is $r$-partite if there is a partition $V_1\cup \cdots\cup V_r$ of $V$ such that each edge $e \in E$ satisfies $|e\cap V_i| = 1$ for each $i$.
Furthermore, we say such a hypergraph is $n$-balanced if $|V_i| = n$ for each $i$.

For each $v \in V$, we let $E_H(v)$ denote the edges containing $v$, $N_H(v)$ denote the set of vertices contained in the edges in $E_H(v)$ apart from $v$ itself, $\deg_H(v) \coloneqq |E_H(v)|$, and $\Delta(H) \coloneqq \max_{u\in V}\deg_H(u)$.
For any $S\subseteq V(H)$, we let $\deg_H(S)$ be the number of edges in $H$ containing $S$, and we define $\Delta_j(H) \coloneqq \max_{|S| = j}\deg_H(S)$ and $\delta_j(H) \coloneqq \min_{|S| = j, \deg_H(S) \geq 1}\deg_H(S)$.

For an $r$-uniform $r$-partite hypergraph $H = (V_1\cup\cdots\cup V_r, E)$, the $r$-partite complement of $H$ is the hypergraph $H^c = (V_1\cup\cdots\cup V_r, E')$, where $E'$ contains all valid edges not in $E$ (an edge is valid if it has exactly one endpoint in each partition $V_i$).

An independent set is a set $I \subseteq V$ containing no edges, and a proper $q$-coloring is a function $\phi\,:\, V \to [q]$ such that $\phi^{-1}(i)$ is an independent set for each $i$.
Finally, recall that a matching in a hypergraph is a set of pairwise disjoint edges.
A perfect matching $M$ in $H$ is a matching such that $v$ is contained in some edge in $M$ for every $v \in V(H)$.

For $s \to \infty$, we use the asymptotic notation $O_s(\cdot)$, $o_s(\cdot)$, etc.
We drop the subscript when $s$ is clear from context.

\section{Introduction}\label{section: intro}

\subsection{Background and main result}\label{subsection: background}
Graph coloring problems have long been central topics in combinatorics, with numerous applications in diverse areas such as theoretical computer science, statistical physics, and network theory. 
One fundamental question in this area is determining the chromatic number of a given graph $G$ (denoted $\chi(G)$), i.e., the minimum number of colors $q$ required for a proper $q$-coloring.
As determining $\chi(G)$ is a known NP-hard problem~\cite{garey1974some}, one aims to study the typical behavior of $\chi(G)$ by considering random structures, such as \ER graphs.
Determining the asymptotic behavior of the chromatic number of such graphs has used several innovative and powerful techniques.

The study of hypergraph coloring extends this line of research to this more general setting. 
Classical results have established bounds on the chromatic number of hypergraphs, particularly in $r$-uniform hypergraphs, under various structural constraints ~\cite{frieze2013coloring,cooper2016coloring,li2022chromatic}.
Meanwhile, the asymptotic behavior of the chromatic number of $r$-uniform \ER hypergraphs was first studied in the dense regime by Shamir \cite{shamir1989chromatic}, and later completely determined in all regimes by Krivelevich and Sudakov \cite{krivelevich1998chromatic}.

When considering $r$-partite hypergraphs, $\bigcup_{i \in J}V_i$ is an independent set for any $J \subsetneq [r]$, and a trivial $2$-coloring exists.
Namely, color the vertices in $V_1$ blue, and color all others red.
A natural question to ask is the following: is the problem still trivial when we impose additional constraints on the color classes?
We will consider the constraint of \textit{balancedness}.

\begin{Definition}\label{def:independent}
    Let $H = (V_1 \cup \cdots \cup V_r, E)$ be an $r$-uniform $r$-partite hypergraph.
    An independent set $I$ in $H$ is \textbf{balanced} if $|I\cap V_i| = |I\cap V_j|$ for each $i, j \in [r]$.
    For the largest balanced independent set $I \subseteq V(H)$, we let $\alpha_b(H) \coloneqq |I|$ denote the \textbf{balanced independence number}.

    A proper $q$-coloring $\phi\,:\,V(H) \to [q]$ is \textbf{balanced} if the color classes $\phi^{-1}(i)$ are balanced independent sets for each $i \in [q]$.
    The \textbf{balanced chromatic number} (denoted $\chi_b(H)$) is the minimum number of colors required for a balanced coloring.
\end{Definition}

Balanced independent sets in bipartite graphs were first introduced by Ash in \cite{ash1983two}, where the author aimed to determine a sufficient condition for Hamiltonicity of bipartite graphs.
Feige and Kogan first studied balanced colorings for bipartite graphs \cite{feige2010balanced}.
They showed that any $n$-balanced bipartite graph $G$ with maximum degree $\Delta$ sufficiently large satisfies $\chi_{b}(G) = O\left(\Delta/\log \Delta\right)$ \cite{feige2010balanced}.
Recently, through a surprisingly simple local algorithm, Chakraborti showed that the hidden constant in the $O(\cdot)$ can be taken to be $1 + o_\Delta(1)$ \cite{chakraborti2023extremal}.
This matches the so-called ``shattering threshold'' for coloring triangle-free graphs.
Indeed, a result of Perkins and the second author implies that no local algorithm can obtain a better result \cite{perkins2024hardness}.

Extending the local algorithm of Chakraborty, the first author proved an upper bound on $\chi_b(H)$ for sparse $n$-balanced $r$-uniform $r$-partite hypergraphs for all uniformities (see Theorem~\ref{theorem: balanced coloring deterministic}) \cite{dhawan2023balanced}, recovering the aforementioned result of Chakraborty when $r = 2$.
As a consequence of earlier work of the authors \cite{dhawan2024low}, improving the constant factor in Theorem~\ref{theorem: balanced coloring deterministic} would require new ideas as local algorithms fail to do so.

For random bipartite graphs, it is mentioned in \cite[\S6]{chakraborti2023extremal}, without a proof, that the balanced chromatic number of the \ER bipartite graph $G_{\mathrm{bip}}(n, d/n)$ is concentrated around $d/(2\log d)$ with high probability\footnote{Throughout this paper, ``with high probability'' means with probability at least $1 - 1/\mathsf{poly}(n)$.}.
Our main result concretely establishes this claim and extends beyond to higher uniformities.
Before we state our result, let us define the \ER model for multipartite hypergraphs:

\begin{Definition}\label{definition: models}
    Let $n,\, r \in \N$ such that $n \geq r \geq 2$.
    For $p \in [0, 1]$, we construct the $r$-uniform $r$-partite hypergraph $H \sim \Hrrnp$ on vertex set $[n]\times [r]$ by including each
    \[e \in V_1 \times \cdots \times V_r,\]
    in $E(H)$ independently with probability $p$.
    Here, $V_i = [n] \times \set{i}$.
\end{Definition}

With this definition in hand, we are ready to state our result.

\begin{Theorem}\label{theorem: main result}
    For all $r \geq 2$ and $\eps \in (0, 1)$, there exists $d_0 \in \N$ such that the following holds for all $d \geq d_0$.
    There exists $n_0 \in \N$ such that for any $n \geq n_0$ and $H \sim \Hrrnp$ for $p = d/n^{r-1}$, we have
    \[(1 - \eps)\left(\frac{r-1}{r}\cdot\frac{d}{\log d}\right)^{\frac{1}{r-1}} \,\leq\, \chi_b(H) \,\leq\, (1 + \eps)\left(\frac{r-1}{r}\cdot\frac{d}{\log d}\right)^{\frac{1}{r-1}},\]
    with high probability.
\end{Theorem}

We note that our result matches the asymptotic behavior of the ordinary chromatic number of sparse $r$-uniform \ER hypergraphs \cite{krivelevich1998chromatic}.
At first glance, this may not seem surprising, since the balanced coloring problem seems similar to ordinary coloring.
However, $\chi_b(\cdot)$ behaves quite differently from the ordinary chromatic number, even for $r = 2$ !
We outline a couple of distinctions below:
\begin{itemize}
    \item A simple observation is that a proper balanced coloring exists only when the $r$-uniform $r$-partite hypergraph is $n$-balanced. 
    Furthermore, not all such hypergraphs admit a balanced coloring (e.g., the complete multipartite hypergraph).
    This is in stark contrast to ordinary hypergraph coloring, which always satisfies $\chi(H) = O\left(\Delta(H)^{1/(r-1)}\right)$.

    \item As observed by Feige and Kogan in \cite{feige2010balanced}, removing an independent set will not increase the chromatic number of the remaining graph, while removing a balanced independent set from a bipartite graph might yield a graph with higher balanced chromatic number (or even a graph that is not balanced colorable). Similar arguments extend to the setting of $r \geq 3$.
    
\end{itemize}
This suggests that our result cannot be derived merely by adapting the arguments from \cite{krivelevich1998chromatic} (or \cite{luczak1991chromatic} for $r = 2$) to the balanced setting, which indeed is the case.

The rest of this introduction is structured as follows: in \S\ref{subsection: proof_overview}, we provide an overview of our proof techniques; and in \S\ref{subsection: future}, we outline potential future directions of inquiry.

\subsection{Proof Overview}\label{subsection: proof_overview}

In this section, we outline the core idea of the proof.
As the analysis is quite involved, we will ignore a number of technical details in this informal overview (in particular, the actual arguments and rigorous definitions given in the rest of the paper may slightly diverge from how they are described here). 
At the same time, we will attempt to intuitively explain the key ideas of the proof and especially highlight the differences in our approach and that of \cite{krivelevich1998chromatic, luczak1991chromatic}.
For this section, we fix the following parameters:
\[p = \frac{d}{n^{r-1}}, \qquad H \sim  \Hrrnp, \qquad q = \left(\frac{r - 1}{r}\,\frac{d}{\log d}\right)^{\frac{1}{r-1}}.\]

The asymptotic order of $\alpha_b(H)$ was determined in earlier work of the authors.
In particular, the following can be obtained as a corollary to \cite[Theorem 1.8]{dhawan2024low}:

\begin{Theorem}\label{theorem: stat thresh balanced}
    For all $r \geq 2$ and $\eps \in (0, 1)$, there exists $d_0 \in \N$ such that the following holds for all $d \geq d_0$.
    There exists $n_0 \in \N$ such that for any $n \geq n_0$ and $H \sim \Hrrnp$ for $p = d/n^{r-1}$, we have
    \[(1 - \eps)rn\left(\frac{r}{r-1}\frac{\log d}{d}\right)^{\frac{1}{r-1}} \,\leq\, \alpha_{b}(H) \,\leq\,  (1 + \eps)rn\left(\frac{r}{r-1}\frac{\log d}{d}\right)^{\frac{1}{r-1}},\]
    with high probability.
\end{Theorem}

The lower bound on $\chi_b(H)$ in Theorem~\ref{theorem: main result} follows by the upper bound on $\alpha_b(H)$ above.
Hence, most of our effort lies in proving the upper bound.
We will employ the so-called \textit{expose-and-merge} strategy first introduced by Matula in his work on coloring dense \ER graphs \cite{matula1987expose}.
This technique has been employed in a number of (hyper)graph coloring arguments since (see, e.g., \cite{luczak1991chromatic, krivelevich1998chromatic, mcdiarmid1990chromatic}; see also \cite[Ch. 10]{beineke2015topics} for a textbook treatment of the argument).
Informally, the strategy for hypergraph coloring follows the following steps (for brevity, we do not introduce new notation, but it should be understood that $H$ below is an \ER hypergraph with edge-density $p$):
\begin{enumerate}[label=(\normalfont{}EM\arabic*), labelindent=\parindent, wide]
    \item\label{EM: smm} Show that a $\approx d^{-\eps}$ proportion of $V(H)$ is colorable with $\approx q^{1-\eps}$ colors with high probability.

    \item\label{EM: expose} Conduct the following procedure $\approx q^\eps$ times: select a uniformly random induced subhypergraph $\tilde H$ of $H$ on $\approx n/d^{\eps}$ vertices and apply \ref{EM: smm} to color a large fraction of $V(\tilde H)$; discard the colored vertices and repeat. (Note that every induced subhypergraph of $H$ is itself an instance of an \ER hypergraph.)

    \item\label{EM: final blow} With high probability, the above procedure yields a proper $q$-coloring of all but at most $\approx n/d^\eps$ vertices. Complete the coloring with a fresh set of $o(q)$ colors.
\end{enumerate}
The argument for \ref{EM: smm} follows a delicate application of the second moment method.
The expose-and-merge step is \ref{EM: expose}.
In particular, we repeatedly \textit{expose} the edges of subhypergraphs $\tilde H$, color a large portion of $\tilde H$, and \textit{merge} the colorings.
Let $U \subseteq V(H)$ be the uncolored vertices in \ref{EM: final blow} and let $H' \coloneqq H[U]$.
One can show that $\Delta \coloneqq \Delta(H') = o(d/\log d)$ with high probability.
The proof is completed by a greedy argument for graphs (using an additional $\Delta + 1$ colors), or a Lov\'asz Local Lemma argument for hypergraphs (using an additional $O\left(\Delta^{1/(r-1)}\right)$ colors).

Our proof follows the strategy outlined above.
Lemma~\ref{lemma: almost coloring} provides the formal statement of \ref{EM: smm}, Lemma~\ref{lemma: expose and merge} provides the statement of \ref{EM: expose}, and Lemma~\ref{lemma: color small subsets} provides the statement of \ref{EM: final blow}.
The first two steps (\ref{EM: smm}--\ref{EM: expose}) follow a conventional approach (see, e.g., \cite{krivelevich1998chromatic, luczak1991chromatic}), whereas the final step \ref{EM: final blow}—proving that the induced subhypergraph on the remaining vertices admits a balanced coloring with few colors—is notably more challenging in the balanced setting.

Note that balancedness is a global property.
In particular, it is unclear if there exists a Lov\'asz Local Lemma argument for a greedy bound on $\chi_b(H')$.
The bound $\chi_b(H') = O(r\Delta(H'))$ from \cite[Lemma 2.2]{dhawan2023balanced} is insufficient as the high probability bound on $\Delta(H')$ we can prove is too large for our desired result (even for $r = 2$).
Furthermore, the observation of Feige and Kogan mentioned earlier implies that it is unclear whether or not $H'$ is balanced colorable to begin with.
To summarize, there are two obstacles we need to overcome:
\begin{enumerate}[label=(\normalfont{}Obs\arabic*), labelindent=\parindent, wide]
    \item\label{obstacle1} Is $H'$ balanced colorable?
    \item\label{obstacle2} If so, is $H'$ balanced colorable with few colors?
\end{enumerate}

Let us first discuss \ref{obstacle1}.
As a result of \cite[Fact 2.1]{dhawan2023balanced}, it is enough to show that the multipartite complement of $H'$ admits a perfect matching.
Aharoni, Georgakopoulos, and Spr\"ussel showed that a $k$-balanced $r$-uniform $r$-partite hypergraph $F$ contains a perfect matching if $\delta_{r-1}(F) \geq k/2$ \cite[Theorem 2]{aharoni2009perfect}.
In Claim~\ref{claim: few edges small s}, we establish that $\Delta_{r-1}(H') \leq n'/2$ with high probability, overcoming \ref{obstacle1} ($n' \in \N$ is such that $H'$ is $n'$-balanced).

We can resolve \ref{obstacle2} by leveraging the following theorem, which is of independent interest (see \S\ref{section: FK} for the proof):

\begin{Theorem}\label{theorem: hypergraph version of FK result}
    For any $r \geq 2$, let $H = (V_1\cup \cdots \cup V_r, E)$ be an $n$-balanced $r$-uniform $r$-partite hypergraph satisfying $|E| = dn$.
    If $H$ is balanced colorable, then $\chi_b(H) \leq \max\set{2,\, r(r-1)d + 1}$.
\end{Theorem}

The above result generalizes~\cite[Theorem 2.1]{feige2010balanced}, who considered the case when $r = 2$.
In Claim~\ref{claim: few edges small s}, we show the average degree of $H'$ is sufficiently small with high probability. 
Combining this with Theorem~\ref{theorem: hypergraph version of FK result} completes the proof.
We remark that the actual argument for \ref{EM: final blow} is a bit more intricate, involving an iterative coloring procedure culminating in the above result.

\subsection{Future directions}\label{subsection: future}

We conclude this introduction with a few potential future directions of inquiry.
Our work serves as a stepping stone in this area, and we believe there are many promising questions worthy of exploration.

\paragraph{$\Vec{\boldgamma}$-Balanced Colorings.}
We propose a variant of coloring multipartite hypergraphs that involves generalizing the notion of ``balanced''. 
Perkins and the second author proposed the notion of $\gamma$-balanced independent sets for $n$-balanced bipartite graphs in \cite{perkins2024hardness}, which the authors of this manuscript generalized to $\Vec{\boldgamma}$-balanced independent sets for $n$-balanced $r$-uniform $r$-partite hypergraphs in \cite{dhawan2024low}.
Here, rather than considering an independent set $I$ with an equal proportion of vertices in each partition, the proportion of vertices of $I$ within each partition is specified by a vector $\Vec{\boldgamma} \in (0, 1)^r$ such that $\sum_{i=1}^r \gamma_i=1$ and $|I \cap V_j| = \gamma_j|I|$ for all $j \in [r]$.
We can therefore generalize $\chi_{b}(H)$ to $\chi_{\Vec{\boldgamma}}(H)$, which denotes the minimum number of colors required to properly color $H$ such that each color class is a $\Vec{\boldgamma}$-balanced independent set.
A natural follow-up question to our results is whether we can establish the asymptotic value of $\chi_{\Vec{\boldgamma}}(H)$ for a sparse \ER multipartite hypergraphs.
The lower bound follows from \cite[Theorem 1.8]{dhawan2024low} and so one can focus on the upper bound.
The key challenge lies in proving a greedy bound akin to Theorem~\ref{theorem: hypergraph version of FK result}.

\paragraph{Balanced $\mathbf{\beta}$-Colorings.}
An independent set in a hypergraph is a subset of vertices $I$ that spans no edge.
For $r \geq 3$, one can tighten this condition by restricting the amount of an edge that can be contained in $I$.
In particular, for $\beta \in \set{1, \ldots, r-1}$, we say that $I$ is $\beta$-independent if $|I \cap e| \leq \beta$ for each edge $e \in E(H)$.
For the extreme values of $\beta$ ($\beta = 1$ and $\beta = r-1$), such sets are referred to as \textit{strong} and \textit{weak} independent sets, respectively.
The $\beta$-chromatic number was introduced in \cite{krivelevich1998chromatic}: for a fixed integer $\beta \in \set{1, \ldots, r-1}$, the $\beta$-chromatic number is the minimum number of colors required to properly color the vertices of $H$ such that each color class forms a $\beta$-independent set.
Krivelevich and Sudakov determined the asymptotic behavior of the $\beta$-chromatic number for all $\beta \in \set{1, \ldots, r-1}$ \cite{krivelevich1998chromatic}. 
Since we have already established the weak balanced chromatic number asymptotically, a natural question to ask is whether we can extend this to the \textit{balanced $\beta$-chromatic number}, where every color class is a balanced $\beta$-independent set. 
Additionally, combining with the previous problem, can we determine the asymptotics of the $\Vec{\boldgamma}$-balanced $\beta$-chromatic number?

\paragraph{Two-Point Concentration.}
In the 1970s, Bollob\'as and Erd\H{o}s \cite{bollobas1976cliques} and Matula \cite{matula1970complete} independently showed that $\alpha(G(n, p))$ is concentrated on two values when $p$ is a constant.
More recently, Bohman and Hofstad extended this to nearly all possible values of $p$ \cite{bohman2024two, bohman2024note}.
When considering colorings, an influential series of works starting with Shamir and Spencer \cite{shamir1987sharp} and culminating in the result of Alon and Krivelevich \cite{alon1997concentration} establish two-point concentration of the chromatic number of $G(n, p)$ for a wide range of $p$ (see also \cite{luczak1991note, achlioptas2005two, coja2016chromatic, surya2022concentration}).
We ask whether the same holds for balanced independent sets and colorings of bipartite graphs.
(Note that point concentration is unknown for the ordinary independence and chromatic numbers of hypergraphs, which is a question of independent interest.)

The rest of the paper is structured as follows: in \S\ref{section: prelim}, we will introduce some preliminary facts and tools that will assist with our proofs; in \S\ref{section: FK}, we will prove Theorem~\ref{theorem: hypergraph version of FK result}; in \S\ref{section: proof}, we will prove Theorem~\ref{theorem: main result} modulo a few technical lemmas, which we will prove in \S\S\ref{section: color most}--\ref{section: color small}.

\section{Preliminaries}\label{section: prelim}

In this section, we collect several results that will be useful in our analysis. 
We first present a bound on the balanced chromatic number of sparse $r$-uniform $r$-partite hypergraphs from earlier work of the first author.
We will employ this result in the iterative coloring procedure for \ref{EM: final blow}.

\begin{Theorem}[{\cite[Theorem 1.7]{dhawan2023balanced}}]\label{theorem: balanced coloring deterministic}
    For all $\eps > 0$ and $r \geq 2$, the following holds for $\Delta \in \N$ sufficiently large and some constant $C\coloneqq C(\eps, r) > 0$.
    For $n \geq C\Delta^{5 + 1/(r-1)}\log\Delta$, let $H$ be an $n$-balanced $r$-uniform $r$-partite hypergraph of maximum degree $\Delta$.
    Then, 
    \[\chi_{b}(H) \leq \left(\left(r - 1 + \eps\right)\frac{\Delta}{\log \Delta}\right)^{1/(r-1)}.\]
\end{Theorem}
 
Next, we recall a fundamental fact relating balanced colorability to perfect matchings. 

\begin{Fact}[{\cite[Fact 2.1]{dhawan2023balanced}}]\label{fact:coloring_matching}
    An $r$-uniform $r$-partite hypergraph $H$ has a balanced coloring if and only if its $r$-partite complement has a perfect matching.
\end{Fact}

Additionally, we prove a simple lemma ensuring balanced colorability when $\Delta_{r-1}(H)$ is sufficiently small. 
This will be useful in the arguments for \ref{EM: final blow}.

\begin{lemma}\label{lemma: colorable small r-1 degree}
    If $H$ is an $n$-balanced $r$-uniform $r$-partite hypergraph satisfying $\Delta_{r-1}(H) \leq n/2$, then $H$ is balanced colorable.
\end{lemma}

\begin{proof}
    Consider the $r$-partite complement $H^c$ of $H$.
    Note that $\delta_{r-1}(H^c) \geq n - \Delta_{r-1}(H) \geq n/2$.
    It follows from \cite[Theorem 2]{aharoni2009perfect} that $H^c$ contains a perfect matching $M \coloneqq \set{e_1, \ldots, e_n}$.
    The claim now follows by Fact~\ref{fact:coloring_matching}.
\end{proof}

We will also take advantage of Talagrand's inequality, which will be instrumental in our probabilistic arguments for \ref{EM: smm}.
The original version in \cite[\S10.1]{MolloyReed} contained an error which was rectified in a subsequent paper of the same authors. 
We state the version from \cite{molloy2014colouring} here.

\begin{Theorem}[{Talagrand's Inequality; \cite{molloy2014colouring}}]\label{theo:Talagrand}
    Let $X$ be a non-negative random variable, not identically $0$, which is a function of $m$ independent trials $T_1$, \ldots, $T_m$. Suppose that $X$ satisfies the following for some $\lambda$, $\beta > 0$: 
    \begin{enumerate}[label=(\normalfont{}T\arabic*)]
        \item Changing the outcome of any one trial $T_i$ can change $X$ by at most $\lambda$.
        \item For any $s>0$, if $X \geq s$ then there is a set of at most $\beta s$ trials that certify $X$ is at least $s$.
    \end{enumerate}
    Then for any $\xi \geq 0$, we have
    \begin{align*}
        \P\Big[\big|X-\E[X]\big| \geq \xi + 20\lambda\sqrt{\beta\E[X]} + 64\lambda^2\beta\Big]\leq 4\exp{\left(-\frac{\xi^2}{8\lambda^2\beta(\E[X] + \xi)}\right)}.
    \end{align*}
\end{Theorem}
\section{Proof of Theorem~\ref{theorem: hypergraph version of FK result}}\label{section: FK}

Let us restate the result for convenience.

\begin{Theorem*}[Restatement of Theorem~\ref{theorem: hypergraph version of FK result}]
    For any $r \geq 2$, let $H = (V_1\cup \cdots \cup V_r, E)$ be an $n$-balanced $r$-uniform $r$-partite hypergraph satisfying $|E| = dn$.
    If $H$ is balanced colorable, then $\chi_b(H) \leq \max\set{2,\, r(r-1)d + 1}$.
\end{Theorem*}

The following fact will be useful for our proof.

\begin{Fact}\label{fact: components}
    Let $H$ be an $n$-vertex $r$-uniform hypergraph having $m$ edges.
    Then, $H$ contains at least $n - (r-1)m$ connected components.
\end{Fact}

\begin{proof}
    The claim trivially holds for $m = 0$.
    If $m > 0$, removing a single edge $e$ from $E(H)$ can create at most $r-1$ new connected components.
    The result now follows by induction.
\end{proof}

We begin by showing the result holds for sufficiently sparse graphs.

\begin{lemma}\label{lemma: 2 colorable}
    Let $H = (V_1\cup \cdots \cup V_r, E)$ be an $n$-balanced $r$-uniform $r$-partite hypergraph satisfying $|E| < n$.
    Then $\chi_b(H) \leq 2$.
\end{lemma}

\begin{proof}
    By Fact~\ref{fact: components}, $H$ has at least $rn - (r-1)|E| \geq n + 1$ components.
    Without loss of generality, assume $H$ has $n+1$ components $C_1, \ldots, C_{n+1}$.
    We claim that there exists $I \subseteq [n+1]$ such that the following holds:
    \begin{align}\label{eqn: size of intersection with first two parts}
        \left|\left(\bigcup_{i \in I}C_i\right)\cap \left(V_1 \cup V_2\right)\right| = n.    
    \end{align}
    Let us first show how the above implies the result.
    Consider such a set $I$ and let 
    \[U_1 \coloneqq \left(\bigcup_{i \in I}C_i\right)\cap V_1, \qquad \text{and} \qquad U_2 \coloneqq \left(\bigcup_{i \in I}C_i\right)\cap V_2.\]
    Let $J_1 \subseteq V(H)$ be such that $J_1\cap V_1 = U_1$, $J_1 \cap V_2 = V_2\setminus U_2$, and $J_1 \cap V_i$ consists of an arbitrary set of $|U_1|$ vertices for $3 \leq i \leq r$.
    Let $J_2 = V(H) \setminus J_1$.
    We claim that $J_1$ and $J_2$ are balanced independent sets.
    First, we note that $|J_1 \cap V_i| = |U_1|$ for $i \in [r]$ (for $i = 2$, it follows since $|U_1| + |U_2| = n = |V_2|$).
    Similarly, $|J_2 \cap V_i| = n - |U_1| = |U_2|$ for $i \in [r]$.
    Since the components $C_i$ are disconnected, any edge $e \in E$ intersects $V_1 \cup V_2$ at $U_1 \cup U_2$ or at $(V_1\setminus U_1)\cup (V_2 \setminus U_2)$.
    It follows that $J_1$ and $J_2$ are independent.
    As these sets partition $V(H)$, they form a balanced coloring of $H$.

    Let us now show \eqref{eqn: size of intersection with first two parts}.
    Consider the following subgraphs of $H$:
    \[S_1 = C_1, \qquad S_j = S_{j-1} \cup C_j, \quad 1 < j \leq n + 1.\]
    For each $i \in [n+1]$, let $s_i = |V(S_i) \cap (V_1 \cup V_2)|$.
    Note that $0 < s_1 < s_2 < \cdots < s_n < 2n$.
    If all the values $s_i\mod n$ are distinct, then there exists $i$ such that $s_i \equiv 0\mod n$ implying $s_i = n$. 
    The desired set is then $I \coloneqq [i]$.
    On the other hand, if there exist $j < i$ such that $s_i \equiv s_j \mod n$, then $s_i - s_j = n$ and the desired set is $I \coloneqq [i] \setminus [j]$.
\end{proof}

\begin{lemma}\label{lemma: colorable implies many edges}
    Let $H = (V_1\cup \cdots \cup V_r, E)$ be an $n$-balanced $r$-uniform $r$-partite hypergraph satisfying $\chi_b(H) = q$ where $3 \leq q < \infty$.
    Then $|E| \geq \dfrac{(q-1)n}{r(r-1)}$.
\end{lemma}

\begin{proof}
    Let $\phi\,:\,V(H) \to [q]$ be a proper balanced $q$-coloring of $H$ and let $J_i \coloneqq \phi^{-1}(i)$ be the corresponding color classes.
    For each $i \in [q]$ and $j \in [r]$, let $J_i^{(j)} = J_i\cap V_j$.
    Consider the following quantity:
    \[s = \sum_{1 \leq i < j < k \leq q}|E(H[J_i \cup J_j \cup J_k])|.\]
    By our assumption that $\chi_b(H) = q$, we have the following as a result of Lemma~\ref{lemma: 2 colorable}:
    \[|E(H[J_i \cup J_j \cup J_k])| \,\geq\, |J_i^{(1)}| + |J_j^{(1)}| + |J_k^{(1)}|.\]
    (If not, we could compute a balanced $2$-coloring of $H[J_i \cup J_j \cup J_k]$, which would yield a balanced $(q-1)$-coloring of $H$.)
    With this observation in hand, we have
    \[s \,\geq\, \sum_{1 \leq i < j < k \leq q}(|J_i^{(1)}| + |J_j^{(1)}| + |J_k^{(1)}|) \,=\, \binom{q-1}{2}\sum_{i = 1}^q|J_i^{(1)}| \,=\, \binom{q-1}{2}n.\]
    Let us now count $s$ in a different way.
    Let $e \in E(H)$ be arbitrary.
    Note that there must exist $2$ vertices in $e$ that lie in different color classes.
    It follows that $e$ is counted at most $\binom{r}{2}(q-2)$ times.
    In particular, we conclude the following:
    \[\binom{q-1}{2}n \,\leq\, s \,\leq\, |E|\binom{r}{2}(q-2).\]
    Rearranging the above gives the desired result.
\end{proof}

We are now ready to prove Theorem~\ref{theorem: hypergraph version of FK result}.
If $d < 1$, the claim follows by Lemma~\ref{lemma: 2 colorable}.
Suppose $H$ is not balanced $2$-colorable.
Then, as $3 \leq \chi_b(H) = q < \infty$, we have the following as a result of Lemma~\ref{lemma: colorable implies many edges}:
\[|E| \,\geq\, \frac{(q-1)n}{r(r-1)} \,\implies\, q \,\leq\, r(r-1)d + 1,\]
as desired.

\section{Proof of Theorem~\ref{theorem: main result}}\label{section: proof}

Let us restate the result for convenience.

\begin{Theorem*}[Restatement of Theorem~\ref{theorem: main result}]
    For all $r \geq 2$ and $\eps \in (0, 1)$, there exists $d_0 \in \N$ such that the following holds for all $d \geq d_0$.
    There exists $n_0 \in \N$ such that for any $n \geq n_0$ and $H \sim \Hrrnp$ for $p = d/n^{r-1}$, we have
    \[(1 - \eps)\left(\frac{r-1}{r}\cdot\frac{d}{\log d}\right)^{\frac{1}{r-1}} \,\leq\, \chi_b(H) \,\leq\, (1 + \eps)\left(\frac{r-1}{r}\cdot\frac{d}{\log d}\right)^{\frac{1}{r-1}},\]
    with high probability.
\end{Theorem*}

We will avoid taking floors and ceilings wherever needed (this does not affect the argument, but simplifies a number of computations).
The first lemma shows that with high probability a ``large'' subgraph is balanced colorable with ``few'' colors.
This lemma consists of the delicate second moment argument mentioned in \S\ref{subsection: proof_overview} (see \S\ref{section: color most} for the proof).

\begin{lemma}\label{lemma: almost coloring}
    For all $r \geq 2$ and $\varepsilon \in (0,1)$, there exits a constant $d_0$ such that the following holds for $d \geq d_0$.
    There exists $n_0 \in \N$ such that for $n$, $p$, and $q$ satisfying the following:
    \[n \geq n_0, 
     \quad p = \frac{d}{n^{r-1}}, 
     \quad q = \left((1+\eps)\,\frac{r - 1}{r}\,\frac{d^{1 - \eps}}{\log d}\right)^{\frac{1}{r-1}},\]
    $\Hrrnp$ contains an $s$-balanced subset with $s \geq n d^{-\eps/(r-1)}$ which can be properly balanced colored using at most $q$ colors with probability at least $1 - n^{-3}$.
\end{lemma}

With this lemma in hand, we are able to color most of the graph through the expose-and-merge technique introduced by Matula in \cite{matula1987expose} (see \S\ref{section: expose and merge} for the proof).
For a proper partial coloring $\phi$, we let $\dom(\phi)$ denote the set of colored vertices.

\begin{lemma}\label{lemma: expose and merge}
    For all $r \geq 2$ and $\eps \in (0, 1)$, there exists $d_0 \in \N$ such that the following holds for all $d \geq d_0$.
    There exists $n_0 \in \N$ such that for $n$, $p$, and $q$ satisfying the following:
    \[n \geq n_0, \quad p = \frac{d}{n^{r-1}}, \quad q = \left((1+\eps)\frac{r-1}{r}\cdot \frac{d}{\log d}\right)^{\frac{1}{r-1}},\]
    there is a partial balanced $q$-coloring $\phi$ of $\Hrrnp$ satisfying 
    \[1 - \frac{|\dom(\phi)|}{rn} \in [d^{-\eps/(12r)}/2,\,d^{-\eps/(12r)}]\]
    with probability at least $1 - o(1) - d^{-\eps/24}$.
\end{lemma}

The final key lemma shows that all ``small'' subgraphs are balanced colorable with ``few'' colors (see \S\ref{section: color small} for the proof).

\begin{lemma}\label{lemma: color small subsets}
    For all $r \geq 2$ and $\eps \in (0, 1)$, there exists $d_0 \in \N$ such that the following holds for all $d \geq d_0$.
    There exists $n_0 \in \N$ such that for $n \geq n_0$ and $p = d/n^{r-1}$, the following holds for any $H \sim \Hrrnp$ with probability at least $1 - o(1/n)$.
    For every $(nd^{-\eps})$-balanced subset $U$ of $V(H)$, the subgraph $H[U]$ is balanced $d^{\frac{1-\delta/4}{r-1}}$-colorable, where $\delta = \eps(1-\eps)/(r - 1)$.
\end{lemma}

With these tools in hand, we are ready to prove our main result:

\begin{proof}[Proof of Theorem~\ref{theorem: main result}]
    We may assume without loss of generality that $\eps$ is sufficiently small for all of our arguments.
    The lower bound follows from the upper bound in Theorem~\ref{theorem: stat thresh balanced}.
    
    For the upper bound, we note the following by first applying Lemma~\ref{lemma: expose and merge} with $\eps/4$ and then Lemma~\ref{lemma: color small subsets} with $\eps/(48r)$:
    \begin{align*}
        \P\left[\chi_b(H) \leq \left((1+\eps/2)\frac{r-1}{r}\frac{d}{\log d}\right)^{1/(r-1)}\right] \geq 1 - o(1) - d^{-\eps/96} \geq \frac{1}{2}.
    \end{align*}
    To obtain a high probability bound, we will use an argument of Frieze \cite{frieze1990independence}.
    Let $Y \coloneqq Y(H)$ be the random variable corresponding to the minimum value $s$ such that there is an $s$-balanced subset $S \subseteq [n]\times [r]$ satisfying $\chi_b(H[V\setminus S]) \leq q'\coloneqq \left((1+\eps/2)\frac{r-1}{r}\frac{d}{\log d}\right)^{1/(r-1)}$.
    We will analyze $Y$ by applying the vertex exposure martingale on $H$.
    To this end, we let $Y_i \coloneqq Y(H[[i]\times [r]])$ (note that $Y = Y_n$).
    The random variables $Y_1, \ldots, Y_n$ satisfy the following:
    \begin{itemize}
        \item $Y_i \geq 0$, and
        \item $|Y_{i+1} - Y_i| \leq 1$.
    \end{itemize}
    Therefore, for $\mu = \E[Y]$, we have the following by Azuma's inequality (see, e.g., \cite[Theorem 7.2.1]{AlonSpencer}):
    \[\P[|Y - \E[Y]| \geq 3\sqrt{n\log n}] \leq \frac{1}{n^2}.\]
    Since $\P[Y = 0] \geq 1/2$, we conclude that $\E[Y] \leq 3\sqrt{n\log n}$.
    In particular, with high probability, there exists an $s$-balanced set $S$ such that
    \begin{enumerate}
        \item\label{item:color} $H[V\setminus S]$ is balanced $q'$-colorable, and
        \item $s \leq nd^{-\eps}$.
    \end{enumerate}
    Let $\phi$ be the coloring from \ref{item:color}.
    By potentially uncoloring a balanced subset of $\dom(\phi)$, we may apply Lemma~\ref{lemma: color small subsets} to color a set $S' \supseteq S$ using at most $d^{\frac{1 - \eps/(8(r-1))}{r - 1}} \leq q - q'$ colors with high probability.
    Putting together both colorings completes the proof.
\end{proof}

\section{Proof of Lemma~\ref{lemma: almost coloring}}\label{section: color most}

Here we restate the lemma for convenience.

\begin{lemma*}[Restatement of Lemma~\ref{lemma: almost coloring}]
    For all $r \geq 2$ and $\varepsilon \in (0,1)$, there exits a constant $d_0$ such that the following holds for $d \geq d_0$.
    There exists $n_0 \in \N$ such that for $n$, $p$, and $q$ satisfying the following:
    \[n \geq n_0, 
     \quad p = \frac{d}{n^{r-1}}, 
     \quad q = \left((1+\eps)\,\frac{r - 1}{r}\,\frac{d^{1 - \eps}}{\log d}\right)^{\frac{1}{r-1}},\]
    $\Hrrnp$ contains an $s$-balanced subset with $s \geq n d^{-\eps/(r-1)}$ which can be properly balanced colored using at most $q$ colors with probability at least $1 - n^{-3}$.
\end{lemma*}

For the rest of this section, we will fix $\eps$, $r$, $d$, $n$, and $p$ to be as in the statement of Lemma~\ref{lemma: almost coloring}.
Furthermore, we will assume $\eps$ is sufficiently small for all arguments in this section.
To assist with our proof, we define the random variables $X$ and $X^*$ as follows:
\begin{align*}
    X &\coloneqq \max\set{s\in [n]\,:\, \text{there exists an } s\text{-balanced set in } \Hrrnp \text{ that is balanced } q\text{-colorable}}, \\
    X^* &\coloneqq \min(X,\, 5nd^{-\eps/(r-1)}).
\end{align*}
The following claim shows that $X^*$ is concentrated about its mean.

\begin{Claim}\label{claim: concentration}
    $\P\left[\Big|X^{*} - \E[X^{*}]\Big| \geq nd^{-(1+\eps)\eps/(r-1)}/2\right] \leq 4\exp\left(-\dfrac{nd^{-(1+2\eps)\eps/(r-1)}}{500}\right)$.
\end{Claim}

\begin{proof}
    For each $i \in [n]$ let $Z_i$ denote the outcomes of all edges containing $i$ in the first partition.
    Clearly the variables $Z_i$ are independent.
    Furthermore, changing the outcome of a single $Z_i$ can change the value of $X^*$ by at most $1$ (either $i$ gets counted in $X^*$, or it is no longer counted in $X^*$).
    Finally, in order to certify $X^* \geq s$, it is enough to exhibit the outcomes of $s$ trials $Z_i$ for vertices $i$ that lie in a set determined by $X^*$.
    It follows that $X^*$ satisfies the conditions of \hyperref[theo:Talagrand]{Talagrand's inequality} with $\lambda = \beta = 1$.
    Using the fact that $\E[X^*] \leq 5nd^{-\eps/(r-1)}$, we have
    \begin{align*}
        \P\left[\Big|X^{*} - \E[X^{*}]\Big| \geq nd^{-(1+\eps)\eps/(r-1)}/2\right] &\leq \P\left[\Big|X^{*} - \E[X^{*}]\Big| \geq nd^{-(1+\eps)\eps/(r-1)}/4 + 20\sqrt{\E[X^*]} + 64\right] \\
        &\leq 4\exp\left(-\frac{n^2d^{-2(1+\eps)\eps/(r-1)}}{128(nd^{-(1+\eps)\eps/(r-1)}/4 + \E[X^*])}\right) \\
        &\leq 4\exp\left(-\frac{nd^{-(1+2\eps)\eps/(r-1)}}{500}\right),
    \end{align*}
    as claimed.
\end{proof}

Lemma~\ref{lemma: almost coloring} will now follow if we can show that $\Hrrnp$ contains an $s$-balanced subset which is balanced $q$-colorable for $s \geq n\left(d^{-\eps/(r-1)}+ d^{-(1+\eps)\eps/(r-1)}\right)$ with probability greater than $4\exp\left(-\frac{nd^{-(1+2\eps)\eps/(r-1)}}{500}\right)$.
Indeed, comparing with Claim~\ref{claim: concentration}, we have
\[
    \E[X^{*}] \geq n\left(d^{-\eps/(r-1)}+ d^{-(1+\eps)\eps/(r-1)}\right) - nd^{-(1+\eps)\eps/(r-1)}/2 = n\left(d^{-\eps/(r-1)}+ d^{-(1+\eps)\eps/(r-1)}/2\right).
\]
Therefore, 
\begin{align*}
   \P\left(X^{*}<nd^{-\eps/(r-1)}\right) 
   \leq  
   \P\left(X^{*}-\E[X^{*}]<nd^{-(1+\eps)\eps/(r-1)}/2\right)
   \leq 
  4\exp\left(\frac{nd^{-(1+2\eps)\eps/(r-1)}}{500}\right) <\frac{1}{n^3}.
\end{align*}
To this end, we let $\mathcal{C}$ be the collection of collections of $q$ pairwise disjoint $k_0$-balanced sets in $[n]\times [r]$, where 
\[
    k_0 = n\left((1 - \delta)\,\frac{r}{r-1}\,\frac{\log d}{d}\right)^{1/(r-1)}, \qquad \delta = \frac{\eps}{2}.
\]
Note the following:
\[qk_0 = n((1 + \delta -\eps\delta)d^{-\eps})^{1/(r-1)} \geq n\left(d^{-\eps/(r-1)}+ d^{-(1+\eps)\eps/(r-1)}\right).\]
In particular, it is enough to show that there is some collection of $q$ pairwise disjoint $k_0$-balanced sets that are independent in $\Hrrnp$ with the desired probability.

For each $\mathcal{S} = \set{S_1, \ldots, S_q} \in \mathcal{C}$, let $Y_{\mathcal{S}}$ denote the event that each $S_i$ is independent in $H$ and let $Y = \sum_{\mathcal{S} \in \mathcal{C}}\mathbf{1}\set{Y_{\mathcal{S}}}$.
Then $Y>0$ implies that $X^{*} \geq qk_0 \geq n\left(d^{-\eps/(r-1)}+ d^{-(1+\eps)\eps/(r-1)}\right)$.
Thus it boils down to bound $\P(Y>0)$ from below, which we shall do in the following claim.

\begin{Claim}\label{claim: payley zigmund}
    $\P[Y > 0] \geq 4\exp\left(\dfrac{nd^{-(1+2\eps)\eps/(r-1)}}{500}\right)$.
\end{Claim}

\begin{proof}
    We will employ the Paley–Zygmund inequality: $\P(Y>0) \geq \dfrac{\E[Y]^2}{\E[Y^2]}$.
    It is  easy to see that 
    \begin{align*}
       \E[Y] = \prod_{i = 1}^q {n-(i-1)k_0 \choose k_0}^r(1-p)^{k_0^r}.
    \end{align*}
    Moreover, 
    \begin{align*}
    \E[Y^2] &= \sum_{\mathcal{S}, \mathcal{S}' \in \mathcal{C}} \P[Y_{\mathcal{S}}, Y_{\mathcal{S}'}] \\
    &= \sum_{\mathcal{S} \in \mathcal{C}} \P[Y_{\mathcal{S}}]\sum_{\mathcal{S}' \in \mathcal{C}} \P[Y_{\mathcal{S}'} \mid Y_{\mathcal{S}}] \\
    &= \sum_{\mathcal{S} \in \mathcal{C}} \P[Y_{\mathcal{S}}]\sum_{\mathcal{S}' \in \mathcal{C}} \prod_{i = 1}^q(1 - p)^{k_0^r - \sum_{j = 1}^q\prod_{l = 1}^r|S_i'\cap S_j\cap V_l|} \\
    &= (1 - p)^{qk_0^r}\sum_{\mathcal{S} \in \mathcal{C}} \P[Y_{\mathcal{S}}]\sum_{\mathcal{S}' \in \mathcal{C}}(1 - p)^{-\sum_{i = 1}^q\sum_{j = 1}^q\prod_{l = 1}^r|S_i'\cap S_j\cap V_l|} \\
    &= (1 - p)^{qk_0^r}\sum_{\mathcal{S} \in \mathcal{C}} \P[Y_{\mathcal{S}}]\sum_{\mathcal{S}' \in \mathcal{C}}\prod_{i = 1}^q(1 - p)^{-\sum_{j = 1}^q\prod_{l = 1}^r|S_i'\cap S_j\cap V_l|}.
    \end{align*}
    We will build the inner sum by counting the number of choices for the collection $\mathcal{S}'$.
    As we are interested in an upper bound, we will be lax with the counting.

    Note that for the set $S_i'$, the integers $k_{jl} = |S_i'\cap S_j\cap V_l|$ satisfy the following:
    \[0 \leq k_{jl} \leq k_0, \qquad \sum_{j = 1}^qk_{jl} \leq k_0.\]
    Let $\mathcal{K}$ denote all possible sets of $qr$ integers satisfying the above.
    Once $S_1', \ldots, S_{i-1}'$ are selected, there are at most $n - (i-1)k_0$ vertices to select the vertices in $(S_i' \setminus S_j) \cap V_l$.
    With this in mind, we have
    \begin{align*}
        \E[Y^2] &\leq (1 - p)^{qk_0^r}\sum_{\mathcal{S} \in \mathcal{C}} \P[\mathcal{S}]\prod_{i = 1}^q\sum_{K \in \mathcal{K}}(1 - p)^{-\sum_{j' = 1}^q\prod_{l' = 1}^rk_{j'l'}}\prod_{l = 1}^r\binom{n - (i-1)k_0}{k_0 - \sum_{j = 1}^qk_{jl}}\prod_{j = 1}^q\binom{k_0}{k_{jl}} \\
        &= (1 - p)^{qk_0^r}\E[Y]\prod_{i = 1}^q\sum_{K \in \mathcal{K}}(1 - p)^{-\sum_{j' = 1}^q\prod_{l' = 1}^rk_{j'l'}}\prod_{l = 1}^r\binom{n - (i-1)k_0}{k_0 - \sum_{j = 1}^qk_{jl}}\prod_{j = 1}^q\binom{k_0}{k_{jl}}
    \end{align*}
    By the AM-GM inequality and convexity of $x^r$, we have
    \[\prod_{l' = 1}^r k_{j'l'}\leq 
    \left(\frac{\sum_{l' = 1}^rk_{j'l'}}{r}\right)^r 
    \leq \sum_{l = 1}^r\frac{k_{j'l'}^r}{r}.\]
    Plugging this into the earlier expression, we have
    \begin{align*}
        \E[Y^2] &\leq (1 - p)^{qk_0^r}\E[Y]\prod_{i = 1}^q\sum_{K \in \mathcal{K}}\prod_{l = 1}^r\binom{n - (i-1)k_0}{k - \sum_{j = 1}^qk_{jl}}\prod_{j = 1}^q\binom{k_0}{k_{jl}}(1-p)^{-k_{jl}^r/r}
    \end{align*}
    Rewriting out the sum over $\mathcal{K}$ in a more explicit way gives us
    \begin{align*}
    &~~~~~\sum_{K \in \mathcal{K}}\prod_{l = 1}^r\binom{n - (i-1)k_0}{k - \sum_{j = 1}^qk_{jl}}\prod_{j = 1}^m\binom{k_0}{k_{jl}}(1-p)^{-k_{jl}^r/r} \\
    = &\sum_{\substack{k_{1,1}, \ldots, k_{q+1, 1} \\
                            \sum_{j=1}^{q+1}k_{j, 1} = k_0 \\
                            \cdots \\
                            k_{1, r}, \ldots, k_{q+1, r} \\
                            \sum_{j=1}^{q+1}k_{j, r} = k_0}} \prod_{l=1}^{r}
                            {k_0 \choose k_{1, l}}\cdots {k_0 \choose k_{q, l}}{n-(i-1)k_0  \choose k_{q+1, l}}(1-p)^{-\sum_{j=1}^qk_{j,l}^r/r}.
    \end{align*}
    
    To assist with our computations, we define the following quantity:
    \begin{align}\label{eqn:a_ell}
        a_{l,s_l} = 
        \sum_{\substack{k_{1, l}, \ldots, k_{q, l} \\
        \sum_{j=1}^qk_{j,l} = s_l}} {k_0 \choose k_{1, l}} \cdots {k_0 \choose k_{q, l}} (1-p)^{-\sum_{j=1}^qk_{j,l}^r/r}
    \end{align}
    With this in hand, we may simplify as follows:
    \begin{align}
        \frac{ \E[Y^2]}{ \E[Y]^2} 
        &\leq 
        \prod_{i=1}^q
        \sum_{\substack{k_{1,1}, \ldots, k_{q+1, 1}  \\
                            \sum_{j=1}^{q+1}k_{j, 1} = k_0 \\
                            \cdots \\
                            k_{1, r}, \ldots, k_{q+1, r} \\
                            \sum_{j=1}^{q+1}k_{j, r} = k_0}}
         \prod_{l=1}^{r} \frac{{k_0 \choose k_{1, l}}\cdots {k_0 \choose k_{q, l}}{n-(i-1)k_0  \choose k_{q+1, l}}}{{n-(i-1)k_0  \choose k_{0}} (1-p)^{ \sum_{j=1}^qk_{j,l}^r/r}} \nonumber\\
        &\leq \left(\sum_{\substack{0 \leq s_{1} \leq k_0 \\
                        \ldots \\
                       0 \leq s_{r} \leq k_0
                        }}
         \prod_{l=1}^{r}  a_{l,s_l} \frac{{n-qk_0 \choose k_0 - s_l}}{{n-qk_0 \choose k_0}}\right)^q \nonumber\\
        &\leq
        \left(\sum_{\substack{0 \leq s_{1} \leq k_0 \\
                        \ldots \\
                       0 \leq s_{r} \leq k_0
                        }} \prod_{l=1}^{r} \frac{a_{l,s_l}}{(n-(q+1)k_0)^{s_l}}\left(\frac{k_0!}{(k_0 - s_l)!}\right)\right)^q \nonumber\\
        &=
        \left(\prod_{l=1}^{r}\sum_{0 \leq s_{l} \leq k_0} \frac{a_{l,s_l}}{(n-(q+1)k_0)^{s_l}}\left(\frac{k_0!}{(k_0 - s_l)!}\right)\right)^q \label{eq: a_l s_l}
    \end{align}
    We will now focus on upper bounds for $a_{l, s_l}$.
    To assist with this, we will require a few useful inequalities.
    
    First, suppose $j_1, \ldots, j_t$ are the indices satisfying $k_{j_i, l}$ is at least $f = n\left(\dfrac{\delta\,r\log d}{(r-1)d}\right)^{1/(r-1)}$.
    We have
    \begin{align}\label{eq: size of second part}
        tf \leq \sum_{j}k_{j, l} \leq k_0 \implies t \leq \left(\frac{1-\delta}{\delta}\right)^{1/(r-1)}.
    \end{align}
    
    Next, for every $k_1 \geq k_2 \geq f$ and $k_1 + k_2 \leq s_l \leq k_0$ we claim that
    \[{k_0 \choose k_1}{k_0 \choose k_2}(1-p)^{-\left(k_1\right)^r/r-\left(k_2\right)^r/r} < {k_0 \choose k_1+k_2}(1-p)^{-\left(k_1+k_2\right)^r/r}.\]
    Indeed, when $k_1 + k_2 \leq 0.7k_0$, we have
    \begin{align*}
        &\qquad\frac{{k_0 \choose k_1}{k_0 \choose k_2}(1-p)^{-\left(k_1\right)^r/r-\left(k_2\right)^r/r}}{{k_0 \choose k_1+k_2}(1-p)^{-\left(k_1+k_2\right)^r/r}} \\
        &\leq \frac{k_0!(k_0 - k_1 - k_2)!}{(k_0 - k_2)!(k_0 - k_1)!}{k_1+k_2 \choose k_1}
        \exp\left(-\frac{p}{r}\left(\left(k_1 + k_2\right)^r - \left(k_1\right)^r-\left(k_2\right)^r\right)\right) \\
        &\leq \left(\frac{k_0 - k_1}{k_0 - k_1 - k_2} \, \frac{e(k_1 + k_2)}{k_1}\,\exp\left(-pk_2^{r-1}\right)\right)^{k_1} \\
        &\leq \left(\frac{7e(k_0 - k_1)}{3k_1}\,\exp\left(-\frac{\delta\,r\log d}{r-1}\right)\right)^{k_1} \\
        &\leq \left(\frac{7e}{3}\left(\frac{1 - \delta}{\delta}\right)^{1/(r-1)}\,\exp\left(-\frac{\delta\,r\log d}{r-1}\right)\right)^{k_1} < 1,
    \end{align*}
    for $d$ sufficiently large.
    Whereas for $k_1 + k_2 \geq 0.7k_0$, we have
    \begin{align*}
        \frac{{k_0 \choose k_1}{k_0 \choose k_2}(1-p)^{-\left(k_1\right)^r/r-\left(k_2\right)^r/r}}{{k_0 \choose k_1+k_2}(1-p)^{-\left(k_1+k_2\right)^r/r}} &\leq 2^{k_0}\,2^{k_0}\,\exp\left(-\frac{\delta\,rk_1\log d}{r-1}\right) \\
        &\leq \exp\left(\frac{2\log 2}{0.7}(k_1 + k_2) - \frac{\delta\,rk_1\log d}{r-1}\right) < 1,
    \end{align*}
    for $d$ sufficiently large.
    Hence,
    \begin{align}\label{eq: bound 1}
        {k_0 \choose k_{1, l}} \cdots {k_0 \choose k_{q,l}} (1-p)^{ - \sum_{j=1}^q k_{j,l}^r/r} \leq {k_0 \choose s_l}(1-p)^{-s_l^r/r}
        \leq {k_0 \choose s_l}\exp\left(\frac{(1+\delta)ps_l^r}{r}\right).
    \end{align}
    Furthermore, for every choice of $k_{1}^{'}, \ldots k_{t}^{'}$, such that $\max_i k_{i}^{'} = f$ and $\sum_i k_{i}^{'} = s$ one can get the following inequality
    \begin{align}\label{eq: bound 2}
        \sum_{\substack{k_{1}^{'}, \ldots k_{t}^{'}, \max_i k^{'}_{i} = f\\ \sum_i k^{'}_{i} = s}}
        {k_0 \choose k^{'}_1} \cdots {k_0 \choose k^{'}_t} (1-p)^{ - \sum_i \left(k_{i}^{'}\right)^r/r} 
        \leq {tk_0 \choose s}\exp\left(\frac{(1+\delta)p s f^{r-1}}{r}\right),
    \end{align}
    by combining 
    \begin{align*}
        \sum_{i = 1}^t\frac{\left(k_{i}^{'}\right)^r}{r} 
        = \frac{1}{r}\sum_{i = 1}^t\left(k_{i}^{'}\right)^{r-1}k_{i}^{'} 
        \leq \frac{f^{r-1}}{r}\sum_{i = 1}^t k_{i}^{'}  = \frac{sf^{r-1}}{r},
    \end{align*}
    and 
    \begin{align*}
        \sum_{\substack{k_{1}^{'}, \ldots k_{t}^{'}, \max_i k^{'}_{i} = f\\ \sum_i k^{'}_{i} = s}}
        {k_0 \choose k^{'}_1} \cdots {k_0 \choose k^{'}_t} \leq \binom{t k_0}{s}.
    \end{align*}

    Recall the definition of $a_{l, s_l}$ in \eqref{eqn:a_ell}.
    We will divide the sum into two parts.
    The first part covers the case where all the $k_{j, l}$ are at most $f$.
    For this case, we use \eqref{eq: bound 2}.
    The second part covers the case where at least one of the $k_{j, l}$ is at least $f$.
    We denote by $i$ the sum of all such $k_{j, l}$'s and use \eqref{eq: size of second part}, \eqref{eq: bound 1}, and \eqref{eq: bound 2}.
    This yields the following inequality:
    \begin{align}
        a_{l, s_l} &\leq {qk_0 \choose s_l}\exp\left(\frac{(1+\delta)\delta\,s_l\log d}{r-1}\right) \nonumber \\ 
        &\qquad + n^{\left(\frac{1-\delta}{\delta}\right)^{1/(r-1)}}\sum_{i = f}^{s_l}\binom{k_0}{i} \exp\left(\frac{(1+\delta)pi^r}{r}\right) \binom{qk_0}{s_l - i}\exp\left(\frac{(1+\delta)\delta(s_l-i)\log d}{r-1}\right) \label{eq:a_l_s_l}
    \end{align}
    Let $b_{i, l}$ denote the $i$-th term in the sum above for $f \leq i \leq s_l$.
    We will split the analysis of $b_{i, l}$ into cases.

    \begin{enumerate}[label = \textbf{(Case \arabic*)}, wide]
        \item $f \leq i \leq 0.4k_0$. We have
        \begin{align*}
            \frac{b_{i+1,l}}{b_{i,l}} &= \left(\frac{k_0 - i}{i+1}\right)\exp\left(\frac{(1+\delta)p((i+1)^r - i^r)}{r}\right)
            \left(\frac{s_l - i}{qk_0 - s_l + i + 1}\right)\exp\left(-\frac{(1+\delta)\delta\log d}{r-1}\right) \\
            &\leq \left(\frac{s_l + 1}{i + 1} - 1\right)\left(\frac{k_0 - i}{qk_0 - s_l + i + 1}\right)\exp\left((1+\delta)^2pi^{r-1}-\frac{(1+\delta)\delta\log d}{r-1}\right) \\
            &\leq \left(\frac{s_l + 1}{i + 1} - 1\right)\,\frac{1}{q}\,\exp\left((1+\delta)^2pi^{r-1}-\frac{(1+\delta)\delta\log d}{r-1}\right) \\
            &\leq \frac{k_0}{f}\,\frac{1}{q}\,\exp\left(\frac{(1+\delta)^3r0.4^{r-1}\log d}{r-1}-\frac{(1+\delta)\delta\log d}{r-1}\right). \\
            &\leq 2\left(\frac{(r-1)(1-\delta)}{(1+\eps)rd^{1-\eps}\delta}\right)^{1/(r-1)}\exp\left(\frac{(1+\delta)^3r0.4^{r-1}\log d}{r-1}-\frac{(1+\delta)\delta\log d}{r-1}\right) < 1.
        \end{align*}
        It follows that $b_{i, l}$ is maximum for $i = f$.
        We bound this parameter as follows:
        \begin{align*}
            b_{f, l} &= \binom{k_0}{f} \exp\left(\frac{(1+\delta)pf^r}{r}\right) \binom{qk_0}{s_l - f}\exp\left(\frac{(1+\delta)(s_l-f)\delta\log d}{r-1}\right) \\
            &\leq \binom{(q+1)k_0}{s_l}\exp\left(\frac{(1+\delta)s_l\delta\log d}{r-1}\right) \\
            &\leq \left(\frac{3qk_0d^{\frac{(1+\delta)\delta}{r-1}}}{s_l}\right)^{s_l}.
        \end{align*}

        \item $0.4s_l \leq i \leq (0.99)^{1/(r-1)}s_l$.
        Setting $\alpha = i/s_l$, we have the following since $s_l \leq k_0$:
        \begin{align*}
            b_{i, l} &= \binom{k_0}{i} \exp\left(\frac{(1+\delta)pi^r}{r}\right) \binom{qk_0}{s_l - i}\exp\left(\frac{(1+\delta)(s_l-i)\delta\log d}{r-1}\right) \\
            &\leq \left(\left(\frac{ek_0}{\alpha s_l}\right)^{\alpha}\exp\left(\frac{(1+\delta)p\alpha^rs_l^{r-1}}{r}\right)\left(\frac{eqk_0}{(1 - \alpha)s_l}\right)^{1 - \alpha}\exp\left(\frac{(1+\delta)(1-\alpha)\delta\log d}{r-1}\right)\right)^{s_l} \\
            &\leq \left(\frac{k_0 d^{((1 - \eps)(1 - \alpha) + (1-\delta^2)\alpha^r)/(r-1)}d^{\frac{(1+\delta)\delta}{r-1}}}{s_l}\right)^{s_l}.
        \end{align*}
        Note that $(1 - \eps)(1 - \alpha) + (1-\delta^2)\alpha^r$ is increasing in $\alpha$ for $\alpha \geq 0.4$.
        Furthermore, for $\alpha = (0.99)^{1/(r-1)}$, the expression is strictly less than $1 - \eps$.
        It follows that
        \[b_{i, l} \leq \left(\frac{3qk_0d^{\frac{(1+\delta)\delta}{r-1}}}{s_l}\right)^{s_l}.\]

        \item $(0.99)^{1/(r-1)}s_l \leq i \leq s_l$.
        We will further split this case into two subcases.
        \begin{enumerate}[label = \textit{(Subcase 3\alph*)}, wide]
            \item $s_l \leq (0.9)^{1/(r-1)}k_0$.
            Here, we have
            \begin{align*}
                b_{i, l} &= \binom{k_0}{i} \exp\left(\frac{(1+\delta)pi^r}{r}\right) \binom{qk_0}{s_l - i}\exp\left(\frac{(1+\delta)(s_l-i)\delta\log d}{r-1}\right) \\
                &\leq \binom{2k_0}{s_l} \exp\left(\frac{(1+\delta)ps_l^r}{r}\right) \left(\frac{eqk_0}{(1 - (0.99)^{1/(r-1)})s_l}\right)^{(1 - (0.99)^{1/(r-1)})s_l} d^{\frac{0.01(1+\delta)\delta\,s_l}{r-1}} \\
                &\leq \left(\frac{2ek_0}{s_l}\exp\left(\frac{(1+\delta)ps_l^{r-1}}{r}\right) \left(\frac{qk_0}{s_l}\right)^{0.01} d^{0.01(1+\delta)\delta/(r-1)}\right)^{s_l} \\
                &\leq \left(\frac{2ek_0}{s_l}\exp\left(\frac{(1+\delta)p0.9k_0^{r-1}}{r}\right) \left(\frac{q}{(0.9)^{1/(r-1)}}\right)^{0.01} d^{0.01(1+\delta)\delta/(r-1)}\right)^{s_l} \\
                &\leq \left(\frac{6k_0d^{(0.91 - 0.89\delta^2 - 0.01\delta)/(r-1)}}{s_l}\right)^{s_l} \leq \left(\frac{3qk_0d^{\frac{(1+\delta)\delta}{r-1}}}{s_l}\right)^{s_l}.
            \end{align*}

            \item $s_l \geq (0.9)^{1/(r-1)}k_0$.
            Note the following for $(0.99)^{1/(r-1)}s_l \leq i \leq s_l(1 - d^{-0.3/(r-1)})$:
            \begin{align*}
                \frac{b_{i+1, l}}{b_{i, l}} &=\left(\frac{s_l + 1}{i + 1} - 1\right)\left(\frac{k_0 - i}{qk_0 - s_l + i + 1}\right) \\
                &\qquad\exp\left(\frac{(1+\delta)p((i+1)^r - i^r)}{r}\right)\exp\left(-\frac{(1+\delta)\delta\log d}{r-1}\right) \\
                &\geq \left(\frac{s_l + 1}{s_l(1-d^{-0.3/(r-1)}) + 1} - 1\right)\left(\frac{d^{-0.3/(r-1)} k_0}{qk_0 + 1}\right)\\
                &\qquad \exp\left(\frac{(1-\delta^2)0.89r\log d}{r-1}\right)\exp\left(-\frac{(1+\delta)\delta\log d}{r-1}\right) \\
                &\geq \frac{1}{10}\,\frac{d^{-0.6/(r-1)}}{q}d^{(1-\delta^2)0.89r/(r-1)}d^{\frac{(1+\delta)\delta}{r-1}} \gg 1.
            \end{align*}
            Furthermore, for $s_l(1 - d^{-0.3/(r-1)}) \leq i \leq s_l$, we have
            \begin{align*}
                b_{i, l} &\leq \binom{k_0}{s_l(1 - d^{-0.3/(r-1)})} \exp\left(\frac{(1+\delta)ps_l^r}{r}\right) \binom{qk_0}{s_l\,d^{-0.3/(r-1)}}\exp\left(\frac{(1+\delta)s_l\,d^{-0.3/(r-1)}\delta\log d}{r-1}\right) \\
                &\leq \binom{k_0}{s_l(1 - d^{-0.3/(r-1)})} \exp\left(\frac{(1+\delta)ps_l^r}{r}\right) \left(\frac{eqk_0}{s_l\,d^{-0.3/(r-1)}}\right)^{s_l\,d^{-0.3/(r-1)}}d^{(1+\delta)\delta\,s_ld^{-0.3/(r-1)}/(r-1)} \\
                &\leq \binom{k_0}{k_0 - s_l}\,\frac{\binom{k_0}{k_0 - s_l(1 - d^{-0.3/(r-1)})}}{\binom{k_0}{k_0 - s_l}} \exp\left(\frac{(1+\delta)ps_l^r}{r}\right) d^{0.8s_l\,d^{-0.3/(r-1)}/(r-1)} \\
                &\leq \binom{k_0}{k_0 - s_l}\left(d^{0.3/(r-1)}\right)^{s_l\,d^{-0.3/(r-1)}} \exp\left(\frac{(1+\delta)ps_l^r}{r}\right) d^{0.8s_l\,d^{-0.3/(r-1)}/(r-1)} \\
                &\leq \binom{k_0}{k_0 - s_l}\exp\left(\frac{(1+\delta)ps_l^r}{r}\right) d^{2k_0\,d^{-0.3/(r-1)}/(r-1)}.
            \end{align*}
            Therefore, in this case $b_{i, l}$ is at most the final expression above.
        \end{enumerate}
    \end{enumerate}
    Putting together all cases, we note the following:
    \begin{align*}
        \sum_{i = f}^{s_l}b_{i, l} &\leq \mathbf{1}\set{s_l \geq (0.9)^{1/(r-1)}k_0}(1 - (0.99)^{1/(r-1)})s_l\binom{k_0}{k_0 - s_l}\exp\left(\frac{(1+\delta)ps_l^r}{r}\right) d^{\frac{2k_0d^{-0.3/(r-1)}}{r-1}} \\
        &\qquad + s_l \left(\frac{3qk_0d^{\frac{(1+\delta)\delta}{r-1}}}{s_l}\right)^{s_l}
    \end{align*}

    Plugging this into \eqref{eq:a_l_s_l}, we obtain:
    \begin{align*}
        \frac{a_{l, s_l}}{n^{\left(\frac{1-\delta}{\delta}\right)^{1/(r-1)}}} &\leq \mathbf{1}_{s_l \geq (0.9)^{1/(r-1)}k_0}(1 - (0.99)^{1/(r-1)})s_l\binom{k_0}{k_0 - s_l}\exp\left(\frac{(1+\delta)ps_l^r}{r}\right) d^{\frac{2k_0d^{-0.3/(r-1)}}{r-1}} \\
        &\qquad + s_l \left(\frac{3qk_0d^{\frac{(1+\delta)\delta}{r-1}}}{s_l}\right)^{s_l} + \frac{{qk_0 \choose s_l}d^{(1+\delta)\delta\,s_l/(r-1)}}{n^{\left(\frac{1-\delta}{\delta}\right)^{1/(r-1)}}} \\
        &\leq \mathbf{1}_{s_l \geq (0.9)^{1/(r-1)}k_0}(1 - (0.99)^{1/(r-1)})s_l\binom{k_0}{k_0 - s_l}\exp\left(\frac{(1+\delta)ps_l^r}{r}\right) d^{\frac{2k_0d^{-0.3/(r-1)}}{r-1}} \\
        &\qquad + 2s_l \left(\frac{3qk_0d^{\frac{(1+\delta)\delta}{r-1}}}{s_l}\right)^{s_l}.
    \end{align*}
    In particular, we have
    \begin{align*}
        &\qquad\sum_{0 \leq s_{l} \leq k_0} \frac{a_{l,s_l}}{(n-(q+1)k_0)^{s_l}}\,\frac{k_0!}{(k_0 - s_l)!} \\
        &\leq n^{\left(\frac{1-\delta}{\delta}\right)^{1/(r-1)}}\sum_{0 \leq s_l \leq k_0}2s_l \left(\frac{3qk_0d^{\frac{(1+\delta)\delta}{r-1}}}{s_l(n-(q+1)k_0)}\right)^{s_l}\,\frac{k_0!}{(k_0 - s_l)!} \\
        &\qquad + n^{\left(\frac{1-\delta}{\delta}\right)^{1/(r-1)}}d^{2k_0\,d^{-0.3/(r-1)}/(r-1)}\sum_{j = 0}^{(1 - (0.9)^{1/(r-1)})k_0}(k_0 - j)\binom{k_0}{j}\exp\left(\frac{(1+\delta)p(k_0 - j)^r}{r}\right)\\
        &\qquad \qquad \times\frac{k_0!}{j!(n-(q+1)k_0)^{k_0 - j}}.
    \end{align*}
    Let us consider the first sum above.
    We have
    \begin{align*}
        \sum_{0 \leq s_l \leq k_0}2s_l \left(\frac{3qk_0d^{\frac{(1+\delta)\delta}{r-1}}}{s_l(n-(q+1)k_0)}\right)^{s_l}\,\frac{k_0!}{(k_0 - s_l)!} &\leq \sum_{0 \leq s_l \leq k_0}\left(\frac{6qk_0^2d^{\frac{(1+\delta)\delta}{r-1}}}{s_l(n-(q+1)k_0)}\right)^{s_l} \\
        &\leq \sum_{1 \leq s_l \leq k_0}\left(\frac{12qk_0^2d^{\frac{(1+\delta)\delta}{r-1}}}{s_l\,n}\right)^{s_l} \\
        &\leq k_0\exp\left(\frac{12qk_0^2d^{\frac{(1+\delta)\delta}{r-1}}}{n}\right).
    \end{align*}
    For the second sum, we have
    \begin{align*}
        &\qquad \sum_{j = 0}^{(1 - (0.9)^{1/(r-1))k_0}}(k_0 - j)\binom{k_0}{j}\exp\left(\frac{(1+\delta)p(k_0 - j)^r}{r}\right)\frac{k_0!}{j!(n-(q+1)k_0)^{k_0 - j}} \\
        &\leq k_0\exp\left(\frac{(1+\delta)pk_0^r}{r}\right)\sum_{j = 0}^{(1 - (0.9)^{1/(r-1)})k_0}\frac{n^j}{j!}\binom{k_0}{j}\exp\left(\frac{(1+\delta)p((k_0 - j)^r - k_0^r)}{r}\right)\left(\frac{2k_0}{n}\right)^{k_0} \\
        &\leq k_0\left(\exp\left(\frac{(1+\delta)pk_0^{r-1}}{r}\right)\frac{2k_0}{n}\right)^{k_0}\\
        &\qquad \qquad \times \sum_{j = 0}^{(1 - (0.9)^{1/(r-1)})k_0}\left(\frac{e^2k_0\,n}{j^2} \exp\left(-\frac{(1+\delta)p\sum_{i = 0}^{r-1}k_0^{r - 1 - i}(k_0 - j)^{i}}{r}\right)\right)^j.
    \end{align*}
    First, we note that
    \begin{align*}
        k_0\left(\exp\left(\frac{(1+\delta)pk_0^{r-1}}{r}\right)\frac{2k_0}{n}\right)^{k_0} &\leq \left(d^{(1-\delta^2)/(r-1)}\frac{4k_0}{n}\right)^{k_0} \ll 1.
    \end{align*}
    Next, we have
    \begin{align*}
        \exp\left(-\frac{(1+\delta)p\sum_{i = 0}^{r-1}k_0^{r - 1 - i}(k_0 - j)^{i}}{r}\right) &\leq \exp\left(-\frac{(1+\delta)pk_0^{r - 1}\sum_{i = 0}^{r-1}(0.9)^{i/(r-1)}}{r}\right) \\
        &= \exp\left(-\frac{(1-\delta^2)\log d\sum_{i = 0}^{r-1}(0.9)^{i/(r-1)}}{r-1}\right) \\
        &\leq d^{-\frac{(1-\delta^2)(1 + 0.9(r-1))}{r-1}}.
    \end{align*}
    It follows that
    \begin{align*}
        &\qquad d^{2k_0\,d^{-0.3/(r-1)}/(r-1)}\sum_{j = 0}^{(1 - (0.9)^{1/(r-1)})k_0}\left(\frac{e^2k_0\,n}{j^2d^{\frac{(1-\delta^2)(1 + 0.9(r-1))}{r-1}}}\right)^j \\
        &\leq d^{2k_0\,d^{-0.3/(r-1)}/(r-1)}k_0\exp\left(\sqrt{\frac{k_0\,n}{d^{\frac{(1-\delta^2)(1 + 0.9(r-1))}{r-1}}}}\right),
    \end{align*}
    where we use the fact that $\max_x(A/x^2)^x$ is $\exp(\sqrt{A}/e)$.
    Putting everything together, we conclude that
    \begin{align*}
        \sum_{0 \leq s_{l} \leq k_0} \frac{a_{l,s_l}}{(n-(q+1)k_0)^{s_l}}\,\frac{k_0!}{(k_0 - s_l)!} &\leq n^{\left(\frac{1-\delta}{\delta}\right)^{1/(r-1)}}k_0\exp\left(\frac{12qk_0^2d^{\frac{(1+\delta)\delta}{r-1}}}{n}\right) + \\
        &\qquad n^{\left(\frac{1-\delta}{\delta}\right)^{1/(r-1)}}d^{2k_0\,d^{-0.3/(r-1)}/(r-1)}k_0\exp\left(\sqrt{\frac{k_0\,n}{d^{\frac{(1-\delta^2)(1 + 0.9(r-1))}{r-1}}}}\right) \\
        &\leq \exp\left(\frac{15qk_0^2d^{\frac{(1+\delta)\delta}{r-1}}}{n}\right).
    \end{align*}

    Finally, we have the following:
    \begin{align*}
        \P[Y > 0] \geq \frac{\E[Y]^2}{\E[Y^2]} &\geq \exp\left(-\frac{15rq^2k_0^2d^{\frac{(1+\delta)\delta}{r-1}}}{n}\right) \\
        &\geq 4\exp\left(-\dfrac{nd^{-(1+2\eps)\eps/(r-1)}}{500}\right),
    \end{align*}
    as desired.
\end{proof}

\section{Proof of Lemma~\ref{lemma: expose and merge}}\label{section: expose and merge}

Let us restate the lemma for convenience.

\begin{lemma*}[Restatement of Lemma~\ref{lemma: expose and merge}]
    For all $r \geq 2$ and $\eps \in (0, 1)$, there exists $d_0 \in \N$ such that the following holds for all $d \geq d_0$.
    There exists $n_0 \in \N$ such that for $n$, $p$, and $q$ satisfying the following:
    \[n \geq n_0, \quad p = \frac{d}{n^{r-1}}, \quad q = \left((1+\eps)\frac{r-1}{r}\cdot \frac{d}{\log d}\right)^{\frac{1}{r-1}},\]
    there is a partial balanced $q$-coloring $\phi$ of $\Hrrnp$ satisfying 
    \[1 - \frac{|\dom(\phi)|}{rn} \in [d^{-\eps/(12r)}/2,\,d^{-\eps/(12r)}]\]
    with probability at least $1 - o(1) - d^{-\eps/24}$.
\end{lemma*}

We will use the expose-and-merge technique introduced by Matula in \cite{matula1987expose}.
Before we provide the details of the algorithm, we define a few parameters for $\beta = \frac{5\eps}{12(r-1)}$, $\xi = \frac{\eps}{12}$, and $\gamma = \frac{\eps}{12r}$:
\begin{align*}
    \overline{n} &= nd^{-\beta}, \qquad \overline{d} = d^{1 - (r-1)\beta}, \qquad \overline{q} = \left((1+\xi)\frac{r-1}{r}\frac{\overline{d}^{1 - \xi}}{\log \overline{d}}\right)^{1/(r-1)} \\
    T &= d^{\beta}\overline{d}^{\xi/(r-1)}(1 - d^{-\gamma}/2).
\end{align*}
The following are easy to verify as a result of our choices of parameters:
\begin{align*}
    T\overline{n}\overline{d}^{-\xi/(r-1)} = n(1 - d^{-\gamma}/2), \qquad \text{and} \qquad T\overline{q} \leq q.
\end{align*}
Consider the following algorithm:

\begin{breakablealgorithm}
\caption{Expose and Merge}\label{alg:expose and merge}

\begin{algorithmic}[1]
    \State $E \gets \emptyset$, $F_0 \gets \emptyset$, $W_{0, j} \gets \emptyset$ for each $j \in [r]$
    \For{$i = 1\ldots T$} 
        \State For $j \in [r]$, choose $A_{i, j} \subseteq [n] \setminus W_{i-1, j}$ with $|A_{i, j}| = \overline{n}$ u.a.r.
        \State Let $H_i$ be the \ER random $r$-uniform $r$-partitite hypergraph with vertex partitions $A_{i, j}$ and edge set $E_i$ sampled with probability $p$.
        \State\label{step:failure} Choose a family of disjoint balanced independent sets $\set{R_1^i, \ldots, R_{\overline{q}}^i}$ in $H_i$ such that $\sum_{s = 1}^{\overline{q}} |R_s^i| = r\overline{n}\overline{d}^{-\xi/(r-1)}$. If no such family exists, return \textsf{FAIL}.
        \State $E_i' \gets E_i \setminus F_{i-1}$
        \State $E \gets E \cup E_i'$
        \State $F_i \gets F_{i-1} \cup (A_{i, 1} \times \cdots \times A_{i, r})$
        \State $W_{i, j} \gets W_{i - 1, j} \cup \left(\bigcup_{s = 1}^{\overline{q}} R_s^i \cap A_{i, j}\right)$
    \EndFor
    \State $\overline{F} = ([n] \times [r]) \setminus \left(\bigcup_{i = 1}^{T}F_i\right)$
    \State Let $\overline{E}$ be a $p$-random subset of $\overline{F}$
    \State $E \gets E \cup \overline{E}$
    \State \Return $E$, $\set{R_1^1, \ldots, R_{\overline{q}}^{T}}$
\end{algorithmic}
\end{breakablealgorithm}

A few remarks are in order. 
First, note that the probability that $e \in E$ is precisely $p$ for each $e \in [n] \times [r]$ and so we may treat the hypergraph $H = ([n] \times [r], E)$ as an instance of $\Hrrnp$.
Additionally, we may think of each $H_i$ as an instance of $\mathcal{H}(r, \overline{n}, p)$.
The following claim will bound the probability that Algorithm~\ref{alg:expose and merge} fails.

\begin{Claim}\label{claim: algorithm failure}
    Algorithm~\ref{alg:expose and merge} fails with probability at most $\dfrac{T}{\overline{n}^3}$.
\end{Claim}

\begin{proof}
    Note the following:
    \[\overline{n}^{r-1}p = d^{1 - (r-1)\beta} = \overline{d}.\]
    By Lemma~\ref{lemma: almost coloring}, it follows that the probability of failure at Step~\ref{step:failure} is at most $1/\overline{n}^3$.
    The claim now follows by a union bound over all iterations.
\end{proof}

By Claim~\ref{claim: algorithm failure}, with probability at least 
\[1 - T/\overline{n}^3 \,\geq\, 1 - 1/n^2,\]
Algorithm~\ref{alg:expose and merge} finds $T\overline{q} \leq q$ disjoint balanced sets $R_1^1, \ldots, R_{\overline{q}}^{T}$ satisfying
\[\sum_s\sum_i|R_s^i| = rT\overline{n}\overline{d}^{-\xi/(r-1)} = rn(1 - d^{-\gamma}/2).\]
Note that although $R_s^i$ is independent in each $H_i$, it may not be independent in $H$.
Let $X$ be the random variable denoting the number of edges of $H$ contained in $R_s^i$ for some $s$ and $i$.

\begin{Claim}\label{claim: bound EX}
    $E[X] \leq nd^{-r\gamma}/2$.
\end{Claim}

\begin{proof}
    Let $e \in [n] \times [r]$ be such that $e \in E$ and $e$ is contained in $R_s^i$ for some $s$ and $i$.
    This can only happen if $e \in E(H_j)$ for some $j < i$.
    Since for all $i$ and $j$, we select $A_{i, j}$ from a set of at least $n - T\overline{n}\overline{d}^{-\xi}$ vertices, we may conclude that the vertices of $e$ are contained in both $A_{i, 1} \times \cdots \times A_{i, r}$ and $A_{j, 1} \times \cdots \times A_{j, r}$ with probability at most 
    \[\left(\frac{\overline{n}}{n - T\overline{n}\overline{d}^{-\xi/(r-1)}}\right)^{2r} = (2d^{\gamma - \beta})^{2r}.\]
    Note that for any $1 \leq s \leq \overline{q}$, each balanced subset of $A_{j, 1} \times \cdots \times A_{j, r}$ containing $|R_s^j|$ elements is equally likely to be chosen as $R_s^j$ by the symmetry of the \ER model.
    By Theorem~\ref{theorem: stat thresh balanced}, we may assume that $|R_s^i|$ is a $t$-balanced set for some $t \leq \overline{n}\left((1+\eps)\frac{r}{r-1}\frac{\log \overline{d}}{\overline{d}}\right)^{1/(r-1)}$ for all $i, s$.
    Therefore, since $|A_{j, l}| = \overline{n}$ for each $l \in [r]$, the probability the vertices of $e$ are in the same set $R_s^j$ for some $s$ is at most
    \[\overline{q}\left(\frac{\overline{n}\left((1+\eps)\frac{r}{r-1}\frac{\log \overline{d}}{\overline{d}}\right)^{1/(r-1)}}{\overline{n}}\right)^r \leq \frac{r}{r-1}\frac{\log \overline{d}}{\overline{d}}.\]
    Since the probability that $e$ appears in $H_j$ is $p$, we have
    \begin{align*}
        \E[X] &\leq n^r\binom{T}{2}p(2d^{\gamma - \beta})^{2r}\frac{r}{r-1}\frac{\log \overline{d}}{\overline{d}} \\
        &\leq 2^{2r}nd\,T^2d^{2r\gamma - 2r\beta}\,\overline{d}^{-1}\,\log d \\
        &= 2^rnd^{1 + 2\beta + 2r\gamma - 2r\beta}\,\overline{d}^{-1 + 2\xi}\,\log d \\
        &\leq \frac{1}{2}\,nd^{2r\gamma - (r-1)\beta + 2\xi} \\
        &\leq \frac{1}{2}nd^{-r\gamma},
    \end{align*}
    for $d$ sufficiently large.
\end{proof}

From Claim~\ref{claim: bound EX} and Markov's inequality, we have
\[\P[X \geq nd^{-\gamma}/2] \leq d^{-(r-1)\gamma} \leq d^{-\eps/24}.\]
For each $1 \leq i \leq T$ and $1 \leq s \leq \overline{q}$, remove any edges that appear in $H$ from $R_s^i$ to obtain $\tilde{R}_s^i$.
Clearly, each of the $\tilde R_s^i$ are balanced independent sets in $H$.
Furthermore, with probability at least $1 - 1/n^2 - d^{-\eps/24}$, we have
\[\sum_s\sum_i|\tilde R_s^i| = \sum_s\sum_i|R_s^i| - rX \geq n(1 - d^{-\gamma}),\]
completing the proof of Lemma~\ref{lemma: expose and merge}.

\section{Proof of Lemma~\ref{lemma: color small subsets}}\label{section: color small}

Let us restate the lemma for convenience.

\begin{lemma*}[Restatement of Lemma~\ref{lemma: color small subsets}]
    For all $r \geq 2$ and $\eps \in (0, 1)$, there exists $d_0 \in \N$ such that the following holds for all $d \geq d_0$.
    There exists $n_0 \in \N$ such that for $n \geq n_0$ and $p = d/n^{r-1}$, the following holds for any $H \sim \Hrrnp$ with probability at least $1 - o(1/n)$.
    For every $(nd^{-\eps})$-balanced subset $U$ of $V(H)$, the subgraph $H[U]$ is balanced $d^{\frac{1-\delta/4}{r-1}}$-colorable, where $\delta = \eps(1-\eps)/(r - 1)$.
\end{lemma*}

For the rest of this section, we will fix $\eps$, $r$, $d$, $n$, $p$, and $\delta$ to be as in the statement of Lemma~\ref{lemma: color small subsets}.
We will prove the lemma through a sequence of claims.

\begin{Claim}\label{claim: few edges large s}
    For every $n/d \leq s \leq n/d^{\eps}$ and every $s$-balanced subset $U$ of $V(H)$, the subgraph $H[U]$ contains at most $s\,d^{1-\delta}$ edges.
\end{Claim}

\begin{proof}
    Let $t_s = s\,d^{1-\delta}$.
    By a union bound, the probability that there exists an $n/d \leq s \leq nd^{-\eps}$ and an $s$-balanced subset $U$ of $V(H)$ such that $H[U]$ has at least $t_s$ edges is at most the following:
    \begin{align*}
        \sum_{s = n/d}^{nd^{-\eps}}\binom{n}{s}^r\binom{s^r}{t_s}p^{t_s} &\leq \sum_{s = n/d}^{nd^{-\eps}}\binom{n}{s}^r\left(\frac{e\,s^r\,p}{t_s}\right)^{t_s} \\
        &\leq \sum_{s = n/d}^{nd^{-\eps}}\left(\left(\frac{en}{s}\right)^{rs/t_s}\left(\frac{e\,s^r\,p}{t_s}\right)\right)^{t_s} \\
        &= \sum_{s = n/d}^{nd^{-\eps}}\left(\left(\frac{en}{s}\right)^{r/d^{1-\delta}}\left(e\,d^\delta\left(\frac{s}{n}\right)^{r-1}\right)\right)^{t_s} \\
        &\leq \sum_{s = n/d}^{nd^{-\eps}}\left(e^{2}\,d^\delta\left(\frac{s}{n}\right)^{r-1 - r/d^{1-\delta}}\right)^{s\,d^{1-\delta}} \\
        &\leq \sum_{s = n/d}^{nd^{-\eps}}\left(e^{2}\,d^\delta\left(\frac{1}{d^{\eps}}\right)^{r-1 - r/d^{1-\delta}}\right)^{n\,d^{-\delta}} \\
        &= o\left(\frac{1}{n^2}\right),
    \end{align*}
    for $d,\, n$ sufficiently large.
\end{proof}

\begin{Claim}\label{claim: color most for large s}
    For every $n/d \leq s \leq n/d^{\eps}$ and every $s$-balanced subset $U$ of $V(H)$, consider the subgraph $H_U \coloneqq H[U]$.
    There exists a proper balanced coloring of all but an $s\,d^{-\delta/2}$-balanced subgraph of $H_U$ using at most $d^{\frac{1-\delta/2}{r-1}}$ colors.
\end{Claim}

\begin{proof}
    Fix such a set $U$.
    By Claim~\ref{claim: few edges large s}, we know that $H[U]$ has at most $sd^{1-\delta}$ edges with probability at least $1 - o(1/n^2)$.
    For each partition of $U$, there are at most $sd^{-\delta/2}$ vertices with degree at least $d^{1-\delta/2}$.
    Let $U' \subseteq U$ be obtained by removing the $sd^{-\delta/2}$ highest degree vertices from each partition of $U$.
    It follows by Theorem~\ref{theorem: balanced coloring deterministic} that $H[U']$ can be properly balanced colored with at most $d^{\frac{1-\delta/2}{r-1}}$ colors, completing the proof.
\end{proof}

\begin{Claim}\label{claim: few edges small s}
    For every $n/d^{1 - \delta/2} \leq s \leq n/d$ and every $s$-balanced subset $W$ of $V(H)$, the subgraph $H[W]$ satisfies the following with probability at least $1 - o(1/n^2)$:
    \begin{itemize}
        \item $\Delta_{r-1}(H[W]) \leq s/2$, and
        \item $|E(H[W])| \leq s\log d$.
    \end{itemize}
\end{Claim}

\begin{proof}
    Let us first consider the $(r-1)$-degree condition.
    By a union bound, the probability that such a set $W$ exists is at most
    \begin{align*}
        \sum_{s = n/d^{1 - \delta/2} }^{n/d}\binom{n}{s}^r\,r\,s^{r-1}\,p^{s/2} &\leq \sum_{s = n/d^{1 - \delta/2} }^{n/d}r\left(\frac{en}{s}\right)^{rs}\,s^{r-1}\left(\frac{d}{n^{r-1}}\right)^{s/2} \\
        &\leq r\sum_{s = n/d^{1 - \delta/2} }^{n/d}\left(2ed^{1 - \delta/2} \right)^{rn/d}\,\left(\frac{n}{d}\right)^{r-1}\left(\frac{d}{n^{r-1}}\right)^{n/(2d^{1 - \delta/2})} \\
        &= o\left(\frac{1}{n^2}\right)
    \end{align*}
    for $d,\, n$ sufficiently large.

    Next, let us consider the number of edges in $H[W]$.
    Following a similar argument to the proof of Claim~\ref{claim: few edges large s}, the probability such a set $W$ exists is at most
    \begin{align*}
        \sum_{s = n/d^{1 - \delta/2} }^{n/d}\binom{n}{s}^r\binom{s^r}{s\log d}p^{s\log d} &\leq \sum_{s = n/d^{1 - \delta/2} }^{n/d}\binom{n}{s}^r\left(\frac{e\,s^r\,p}{s\log d}\right)^{s\log d} \\
        &\leq \sum_{s = n/d^{1 - \delta/2} }^{n/d}\left(\left(\frac{en}{s}\right)^{r/\log d}\left(\frac{e\,s^r\,p}{s\log d}\right)\right)^{s\log d} \\
        &= \sum_{s = n/d^{1 - \delta/2} }^{n/d}\left(\left(\frac{en}{s}\right)^{r/\log d}\left(\frac{ed}{\log d}\,\left(\frac{s}{n}\right)^{r-1}\right)\right)^{s\log d} \\
        &\leq \sum_{s = n/d^{1 - \delta/2} }^{n/d}\left(\frac{e^2d}{\log d}\,\left(\frac{s}{n}\right)^{r-1 - r/\log d}\right)^{s\log d} \\
        &\leq \sum_{s = n/d^{1 - \delta/2} }^{n/d}\left(\frac{e^2d}{\log d}\,\left(\frac{1}{d}\right)^{r-1 - r/\log d}\right)^{n\log d/d^{1 - \delta/2} } \\
        &= o\left(\frac{1}{n^2}\right),
    \end{align*}
    for $d,\, n$ sufficiently large.
    The claim now follows by putting together both bounds.
\end{proof}

\begin{Claim}\label{claim: color small s}
    For every $n/d^{1 - \delta/2} \leq s \leq n/d$ and every $s$-balanced subset $W$ of $V(H)$, the subgraph $H[W]$ can be properly balanced colored with at most $r^2\log d$ colors with probability at least $1 - o(1/n^2)$.
\end{Claim}

\begin{proof}
    By Claim~\ref{claim: few edges small s}, with probability at least $1 - o(1/n)$, the subgraph $H[W]$ has at most $s\log d$ edges and maximum $(r-1)$-degree at most $s/2$.
    The claim now follows by Lemma~\ref{lemma: colorable small r-1 degree} and Theorem~\ref{theorem: hypergraph version of FK result}.
\end{proof}

We are now ready to prove Lemma~\ref{lemma: color small subsets}.
Begin with such a set $U$.
Repeatedly apply Claim~\ref{claim: color most for large s}, each time using a fresh set of $d^{\frac{1-\delta/2}{r-1}}$ colors and decreasing the size of the set by a factor of $d^{\delta/2}$.
After at most $(1-\eps)\log d$ iterations, we are left with an $s$-balanced set $W$ such that $n/d^{1-\delta/2} \leq s \leq n/d$.
We may now finish the coloring by applying Claim~\ref{claim: color small s}.
This procedure succeeds with probability at least 
\[1 - o\left(\frac{\log d}{n^2}\right) - o\left(\frac{1}{n}\right) \,=\, 1 - o\left(\frac{1}{n}\right).\]
As there are at most $\Theta(d^{r\eps})$ such sets $U$, the claim holds for all such $U$ with high probability for sufficiently large $n$.
Furthermore, the total number of colors used is at most:
\[(1-\eps)d^{\frac{1-\delta/2}{r-1}}\log d + r^2\log d \,\leq\, d^{\frac{1-\delta/4}{r-1}},\]
for sufficiently large $n,\,d$.

\printbibliography

@BOOK{MolloyReed,
    AUTHOR = " M. Molloy and B. Reed",
    TITLE = "{Graph Colouring and the Probabilistic Method}",
    PUBLISHER = "Springer",
    YEAR = "2002",
}

@book{AlonSpencer,
	author = {N. Alon and J.H. Spencer},
	title = {The Probabilistic Method},
	date = {2000},
	edition = {2},
	publisher = {John Wiley {\&} Sons},
}

@inproceedings{garey1974some,
  title={Some simplified NP-complete problems},
  author={Garey, Michael R and Johnson, David S and Stockmeyer, Larry},
  booktitle={Proceedings of the sixth annual ACM symposium on Theory of computing},
  pages={47--63},
  year={1974}
}

@article{ash1983two,
  title={Two sufficient conditions for the existence of hamiltonian cycles in bipartite graphs},
  author={P. Ash},
  journal={Ars Comb. A},
  volume={16},
  pages={33--37},
  year={1983}
}

@article{shamir1989chromatic,
  title={Chromatic Numbers of Random Hypergraphs and Associated Graphs.},
  author={Shamir, Eli},
  journal={Adv. Comput. Res.},
  volume={5},
  pages={127--142},
  year={1989}
}

@article{krivelevich1998chromatic,
  title={The chromatic numbers of random hypergraphs},
  author={Krivelevich, Michael and Sudakov, Benny},
  journal={Random Structures \& Algorithms},
  volume={12},
  number={4},
  pages={381--403},
  year={1998},
  publisher={Wiley Online Library}
}

@article{aharoni2009perfect,
  title={Perfect matchings in r-partite r-graphs},
  author={R. Aharoni and A. Georgakopoulos and P. Spr{\"u}ssel},
  journal={European Journal of Combinatorics},
  volume={30},
  number={1},
  pages={39--42},
  year={2009},
  publisher={Elsevier}
}

@article{feige2010balanced,
  title={Balanced coloring of bipartite graphs},
  author={U. Feige and S. Kogan},
  journal={Journal of Graph Theory},
  volume={64},
  number={4},
  pages={277--291},
  year={2010},
  publisher={Wiley Online Library}
}

@article{frieze2013coloring,
  title={Coloring simple hypergraphs},
  author={A. Frieze and D. Mubayi},
  journal={Journal of Combinatorial Theory, Series B},
  volume={103},
  number={6},
  pages={767--794},
  year={2013},
  publisher={Elsevier}
}

@article{luczak1991chromatic,
  title={The chromatic number of random graphs},
  author={{\L}uczak, Tomasz},
  journal={Combinatorica},
  volume={11},
  number={1},
  pages={45--54},
  year={1991},
  publisher={Springer}
}

@article{mcdiarmid1990chromatic,
  title={On the chromatic number of random graphs},
  author={McDiarmid, Colin},
  journal={Random Structures \& Algorithms},
  volume={1},
  number={4},
  pages={435--442},
  year={1990},
  publisher={Wiley Online Library}
}

@article{molloy2014colouring,
  title={Colouring graphs when the number of colours is almost the maximum degree},
  author={M. Molloy and B. Reed},
  journal={Journal of Combinatorial Theory, Series B},
  volume={109},
  pages={134--195},
  year={2014},
  publisher={Elsevier}
}

@article{cooper2016coloring,
  title={Coloring sparse hypergraphs},
  author={J. Cooper and D. Mubayi},
  journal={SIAM Journal on Discrete Mathematics},
  volume={30},
  number={2},
  pages={1165--1180},
  year={2016},
  publisher={SIAM}
}

@unpublished{li2022chromatic,
  title={The chromatic number of triangle-free hypergraphs},
  author={Li, Lina and Postle, Luke},
  howpublished={\url{https://arxiv.org/abs/2202.02839} (arXiv preprint)},
  year={2022}
}

@article{dhawan2023balanced,
    author = {Dhawan, Abhishek},
    title = {Balanced Independent Sets and Colorings of Hypergraphs},
    journal = {Journal of Graph Theory},
    volume={109},
    number={1},
    pages={43--51},
    year={2025}
}

@article{chakraborti2023extremal,
  title={Extremal bipartite independence number and balanced coloring},
  author={D. Chakraborti},
  journal={European Journal of Combinatorics},
  volume={113},
  pages={103750},
  year={2023},
  publisher={Elsevier}
}

@inproceedings{perkins2024hardness,
  title={On the hardness of finding balanced independent sets in random bipartite graphs},
  author={Perkins, Will and Wang, Yuzhou},
  booktitle={Proceedings of the 2024 Annual ACM-SIAM Symposium on Discrete Algorithms (SODA)},
  pages={2376--2397},
  year={2024},
  organization={SIAM}
}

@unpublished{dhawan2024low,
  title={The Low-Degree Hardness of Finding Large Independent Sets in Sparse Random Hypergraphs},
  author={Dhawan, Abhishek and Wang, Yuzhou},
  howpublished = {\url{https://arxiv.org/abs/2404.03842} (arXiv preprint)},
  year={2024}
}

@article{frieze1990independence,
  title={On the independence number of random graphs},
  author={Frieze, Alan M},
  journal={Discrete Mathematics},
  volume={81},
  number={2},
  pages={171--175},
  year={1990},
  publisher={Elsevier}
}

@article{matula1987expose,
  title={Expose-and-merge exploration and the chromatic number of a random graph},
  author={Matula, David W},
  journal={Combinatorica},
  volume={7},
  number={3},
  pages={275--284},
  year={1987},
  publisher={Springer}
}

@inproceedings{bohman2024two,
  title={Two-point concentration of the independence number of the random graph},
  author={Bohman, Tom and Hofstad, Jakob},
  booktitle={Forum of Mathematics, Sigma},
  volume={12},
  pages={e24},
  year={2024},
  organization={Cambridge University Press}
}

@unpublished{bohman2024note,
  title={A note on Two-Point Concentration of the Independence Number of $G_{n, m}$},
  author={Bohman, Tom and Hofstad, Jakob},
  howpublished={\url{https://arxiv.org/abs/2410.05420} (arXiv preprint)},
  year={2024}
}

@inproceedings{bollobas1976cliques,
  title={Cliques in random graphs},
  author={Bollob{\'a}s, B{\'e}la and Erd{\H{o}}s, Paul},
  booktitle={Mathematical Proceedings of the Cambridge Philosophical Society},
  volume={80},
  number={3},
  pages={419--427},
  year={1976},
  organization={Cambridge University Press}
}

@inproceedings{matula1970complete,
author={David Matula},
booktitle = {Proceedings of the Second Chapel Hill Conference on Combinatory Mathematics and its Applications},
pages={356–-369},
title = {On complete subgraphs of a random graph},
year = {1970},
}

@article{shamir1987sharp,
  title={Sharp concentration of the chromatic number on random graphs G n, p},
  author={Shamir, Eli and Spencer, Joel},
  journal={Combinatorica},
  volume={7},
  pages={121--129},
  year={1987},
  publisher={Springer}
}

@article{alon1997concentration,
  title={The concentration of the chromatic number of random graphs},
  author={Alon, Noga and Krivelevich, Michael},
  journal={Combinatorica},
  volume={17},
  number={3},
  pages={303--313},
  year={1997},
  publisher={Springer}
}

@article{luczak1991note,
  title={A note on the sharp concentration of the chromatic number of random graphs},
  author={Luczak, Tomasz},
  journal={Combinatorica},
  volume={11},
  number={3},
  pages={295--297},
  year={1991},
  publisher={Springer}
}

@article{achlioptas2005two,
  title={The two possible values of the chromatic number of a random graph},
  author={Achlioptas, Dimitris and Naor, Assaf},
  journal={Annals of Mathematics},
  volume={162},
  pages={1335--1351},
  year={2005}
}

@article{coja2016chromatic,
  title={The chromatic number of random graphs for most average degrees},
  author={Coja-Oghlan, Amin and Vilenchik, Dan},
  journal={International Mathematics Research Notices},
  volume={2016},
  number={19},
  pages={5801--5859},
  year={2016},
  publisher={Oxford University Press}
}

@unpublished{surya2022concentration,
  title={On the concentration of the chromatic number of random graphs},
  author={Surya, Erlang and Warnke, Lutz},
  howpublished={\url{https://arxiv.org/abs/2201.00906} (arXiv preprint)},
  year={2022}
}

@book{beineke2015topics,
  title={Topics in chromatic graph theory},
  author={Beineke, Lowell W and Wilson, Robin J},
  volume={156},
  year={2015},
  publisher={Cambridge University Press}
}

\end{document}